\documentclass[journal,11pt,draftcls,onecolumn,peerreview]{IEEEtran}

\usepackage{url}
\usepackage{amsmath,amssymb,amsthm,bm}
\usepackage{subfig}
\usepackage{cite}
\usepackage[latin1]{inputenc}

\ifCLASSINFOpdf
  \usepackage[pdftex]{graphicx}
  \graphicspath{{figs/pdf/}{figs/jpg/}}
  \DeclareGraphicsExtensions{.pdf,.jpeg,.jpg,.png}
\else
  \usepackage[dvips]{graphicx}
  \graphicspath{{figs/eps/}}
  \DeclareGraphicsExtensions{.eps}
\fi
\newcommand{\abs}[1]{\lvert#1\rvert}
\newcommand{\lemma}{Lemma}
\newtheorem{thm}{\lemma}
\newtheorem{prop}{Property}

\begin{document}

\title{Approximate Maximum Likelihood Source Localization from Range Measurements Through Convex Relaxation}

\author{P\i nar~O\u{g}uz-Ekim,~\IEEEmembership{Student Member,~IEEE,}
        Jo\~{a}o~Gomes,~\IEEEmembership{Member,~IEEE,}
        Jo\~{a}o~Xavier,~\IEEEmembership{Member,~IEEE,}\\
        Marko Sto\v{s}i\'{c},
        and~Paulo~Oliveira,~\IEEEmembership{ Senior Member,~IEEE}
\thanks{The authors are with the Institute for Systems and Robotics -- Instituto Superior T\'ecnico (ISR/IST), Lisbon, Portugal. e-mails: \{poguz,jpg,jxavier,mstosic,pjcro\}@isr.ist.utl.pt.}
\thanks{This research was partially supported by Funda\c{c}\~{a}o para a Ci\^{e}ncia e a Tecnologia (FCT) through ISR/IST plurianual funding with PIDDAC program funds, projects PTDC/EEA-TEL/71263/2006, PTDC/EEA-CRO/104243/2008, CMU-PT/SIA/0026/2009, and grant SFRH/BD/44771/2008.}}

\ifCLASSOPTIONpeerreview
\else
\markboth{IEEE TRANSACTIONS ON SIGNAL PROCESSING}
{O\u{g}uz-Ekim \MakeLowercase{\textit{et al.}}: ML Source Localization Through Convex Relaxation}
\fi

\maketitle

\begin{abstract}
This work considers the problem of locating a single source from noisy range measurements to a set of nodes in a wireless sensor network. We propose two new techniques that we designate as Source Localization with Nuclear Norm (SLNN) and Source Localization with $\ell_1$-norm (SL-$\ell_1$), which extend to arbitrary real dimensions, including 3D, our prior work on 2D source localization formulated in the complex plane. Broadly, our approach is based on formulating a Maximum-Likelihood (ML) estimation problem for the source position, and then using convex relaxation techniques to obtain a semidefinite program (SDP) that can be globally and efficiently solved. SLNN directly approximates the Gaussian ML solution, and the relaxation is shown to be tighter than in other methods in the same class. We present an analysis of the convexity properties of the constraint set for the 2D complex version of SLNN (SLCP) to justify the observed tightness of the relaxation. In terms of global accuracy of localization, SLNN outperforms state-of-the-art optimization-based methods with either iterative or closed-form formulations.
We propose the SL-$\ell_1$ algorithm to address the Laplacian noise case, which models the presence of outliers in range measurements. We overcome the nondifferentiability of the Laplacian likelihood function by rewriting the ML problem as an exact weighted version of the Gaussian case, and compare two solution strategies. One of them is iterative, based on block coordinate descent, and uses SLNN as a subprocessing block. The other, attaining only slightly worse performance, is noniterative and based on an SDP relaxation of the weighted ML problem.
\end{abstract}

\ifCLASSOPTIONpeerreview

\else
\fi
%


\section{Introduction}\label{sec:intro}


Locating a source from range measurements to a set of known reference points (anchors) is a classic problem in many engineering applications (e.g., radar, sonar, GPS), and has received a great deal of attention over the years. Recently, source localization from range measurements has been intensively examined in the context of wireless sensor networks (WSN), where ranges estimated from times of arrival, or from surrogates such as received signal strength, are somewhat unreliable due to the complexity of many WSN propagation environments (e.g., indoor settings with few unobstructed line-of-sight paths).

Spatial information \emph{per se}, or as georeference to other sensor measurements, is crucial in WSN applications and warrants investigation into suitable localization algorithms. While many approaches to source localization based on classical triangulation or heuristic criteria can be found in the WSN literature \cite{Dargie2010, He2003}, our primary focus is on optimization-based methods formally derived from the likelihood function of observations, or related cost functions \cite{Ekim2010, Stoica2008, Beck2008, Xu2010, Chan2006, Cheung2004}. By doing so, we expect to take advantage of the optimality properties of maximum likelihood (ML) estimates to improve the robustness to perturbations in range measurements. We do not consider alternative/complementary measurements such as angles of arrival or time differences of arrival. We also assume cooperative localization scenarios where absolute ranges, as opposed to range differences to a reference sensor, are measured. These can be obtained either by synchronizing clocks and transmitting waveforms from the source at known times (beacon mode), or by initiating the transmission at a reference sensor and measuring the round trip time to the source and back (transponder mode).

Centralized ML algorithms for range-based source localization, which require the transmission of the full data set to a fusion node for processing, are proposed in \cite{Cheung2004, Ekim2010} under Gaussian noise and in \cite{Ekim2011} under Laplacian noise. These resort to semidefinite relaxation (SDR) to alleviate the problem of algorithmic convergence to undesirable local maxima of the likelihood function. A related alternative approach proposed in \cite{Stoica2008} solves a constrained least-squares (LS) problem using squared range (SR) measurements, subject to a quadratic constraint. This was shown to outperform, on average, the ML SDR approach of \cite{Cheung2004} whose relaxed solutions sometimes fail to produce meaningful source position vectors (rank one solutions). Another approach, proposed in \cite{Xu2010}, approximates the ML solution via second-order cone programming and a low-dimensional search.

Distributed algorithms for wireless sensor nodes, where the source location is iteratively determined through in-network processing at individual nodes and communication between neighbours, are also being very actively pursued \cite{Hero2005, So2009, Dimakis2010}. These techniques, however, are not the focus of our work. We also note that source localization can be viewed as a special instance of sensor network localization, where the positions of several sources/sensors are simultaneously determined from pairwise range measurements. Related algorithms based on semidefinite programming (SDP) have been developed for this class of problems \cite{Ekim2011, Ye2004}, and are relevant when there is significant uncertainty in anchor positions (see, e.g., \cite{Yang2010} for a similar SDP approach to source localization with anchor uncertainty using range differences).

This paper develops an alternative to the source localization ML SDR method of \cite{Cheung2004}. We term this approach, originally proposed in \cite{Ekim2010} for 2D localization under Gaussian noise, Source Localization in the Complex Plane (SLCP). Our relaxation for the nonconvex and nonsmooth likelihood function is tighter than the one presented in \cite{Cheung2004}, in the sense that the relaxed solution will more often have (near) rank-1, as required to obtain target coordinates by factorization. SLCP also outperforms the SR-LS method of \cite{Stoica2008}, which iteratively solves a generalized trust-region subproblem and dispenses with factorization of rank-1 matrices, but undergoes some degradation with noisy measurements due to squaring of ranges in the cost function. The degradation of SR-LS becomes more severe in the presence of outliers \cite{Ekim2011}, which commonly affect practical range measurement systems, e.g., when non-line-of-sight propagation occurs.

This paper expands upon the results of \cite{Ekim2010} in several ways:
\begin{enumerate}
\item{We extend the framework of SLCP from 2D localization, which relied on a formulation where target and anchor coordinates were represented as complex numbers, to arbitrary (real) dimensions. We term the new SDR method Source Localization with Nuclear Norm (SLNN), as this norm arises naturally in the cost function of our relaxed optimization problem. Similarly to SLCP, SLNN offers a tight relaxation in most problem instances, and retains a performance advantage over SR-LS.}

\item{We provide a more complete analysis of the accuracy properties of SLCP, whose success in providing tight relaxations relies on certain parametrically defined sets in $\mathbb{R}^2$ being nearly convex. We discuss the convexity of the sets and how to trace the convex hull for any of them, from which convexity can be empirically assessed. For three-anchor scenarios we also examine a search-based alternative to SVD decomposition to extract the source coordinates from the solution of SLCP (a positive semidefinite matrix with near rank-1).}


\item{In \cite{Ekim2011} a modification of SLCP, termed SL-$\ell_1$, was introduced for ML source localization under Laplacian noise. This makes the algorithm robust to outlier measurements, a property that was observed in simulation even for non-Laplacian range errors. In this paper we provide a conceptually similar extension for source localization beyond 2D, consisting of a reformulation of the nondifferentiable log-likelihood function for Laplacian noise as a reweighted version of the Gaussian log-likelihood. We propose both single convex formulations and a simpler iterative optimization algorithm, which repeatedly solves weighted SLNN problems followed by weight refinement.}
\end{enumerate}

The organization of the paper is as follows. In Section~\ref{sec:formulation}, we formulate the ML source location problem under Gaussian or Laplacian noise. In Section \ref{sec:slcp} we derive the SLCP algorithm for 2D localization, we analyze the geometry of the associated optimization problem and the tightness of the relaxation (Subsection \ref{sec:convexhull}), and we propose criteria for factorizing the SDR solution to recover the source coordinates (Subsection \ref{sec:rank1}). In Section \ref{sec:slnn} we derive the SLNN algorithm, which extends SLCP to 3 (and higher) dimensions, and propose an iterative version of this algorithm that can handle Laplacian noise (Section \ref{sec:SLl1}). Section~\ref{sec:simu} illustrates the performances of the algorithms in simulation. Finally, conclusions are drawn in Section~\ref{sec:conclusion}.

Throughout, both scalars and individual position vectors are represented by lowercase letters. Other vectors and matrices are denoted by boldface lowercase and uppercase letters, respectively. Individual components of matrix $\mathbf{X}$ are written as $x_{ij}$ and those of vector $\mathbf{x}$ as $x_{i}$ (the same notation would be used for a hypothetical position vector $x_{i}$, but the distinction between both should be clear from context). The superscript \emph{T} (\emph{H}) denotes the transpose (hermitian) of the given real (complex) vector or matrix, $\langle \cdot, \cdot \rangle$ denotes the inner product of two vectors, and $\text{tr}(\cdot)$ denotes the trace of a matrix. For symmetric matrix $\mathbf{X}$, $\mathbf{X}\succeq 0$ means that $\mathbf{X}$ is positive semidefinite. We denote the Frobenius norm of matrix $\mathbf{X}$ as $\|\mathbf{X}\|_{F} = \sqrt{\text{tr}(\mathbf{X}^{H}\mathbf{X})}$ and its nuclear norm as $\|\mathbf{X}\|_{N} = \text{tr}\bigl( (\mathbf{X}^{H}\mathbf{X})^{\frac{1}{2}} \bigr)$. Below, $\mathbf{I}_m$ is the $m \times m$ identity matrix and $\mathbf{1}_m$ is the vector of $m$ ones. The convex hull of set ${\cal S}$ is denoted by $\text{co}({\cal S})$.


\section {Problem formulation}\label{sec:formulation}

Let $x\in \mathbb{R}^n$ be the unknown source position, $a_i\in \mathbb{R}^n$, $i = 1,..,m$ be known sensor positions (anchors), and $r_{i} = \|x-a_i\|+w_{i}$ be the measured range between the  source and the $i$-th anchor, where $w_{i}$ denotes a noise term with standard deviation $\sigma$. Under i.i.d.\ Gaussian or Laplacian noise maximizing the likelihood of observations for the source localization problem is equivalent to
\begin{equation}\label{eq:costfunction1_2}
\begin{array}{cl}
\text{minimize}& \sum_{i=1}^{m} |\|x-a_i\|^p-r_{i}^p|^q.\\
x
\end{array}
\end{equation}
We will derive the SLCP/SLNN algorithms to (approximately) solve \eqref{eq:costfunction1_2} under Gaussian noise ($p = 1$, $q = 2$), whereas SL-$\ell_1$ will solve it under Laplacian noise ($p = 1$, $q = 1$). The case  ($p = 2$, $q = 2$) is also of interest and corresponds to the cost function used in the SR-LS algorithm of \cite{Stoica2008}, which is used to benchmark our algorithms. Note that the cost function for SR-LS is not a likelihood function, and it arises out of mathematical convenience.

The main difficulties of solving \eqref{eq:costfunction1_2} lie in the fact that this cost function is, in general, nonconvex and multimodal. For $q=1$ it is also nondifferentiable, which poses additional challenges. We address the nonconvexity and multimodality of the cost function in Sections \ref{sec:slcp} and \ref{sec:slnn} by developing convex relaxations that turn out to be tight in most problem instances, thus providing a very good approximation to the true source location. If necessary, the source coordinates can be further refined by iteratively minimizing \eqref{eq:costfunction1_2} starting from the relaxed solution. Refer to \cite{Ekim2011} for one such iterative refinement approach based on the Majorization-Minimization algorithm. We address the nondifferentiability of \eqref{eq:costfunction1_2} for $q=1$ in Section \ref{sec:SLl1} by rewriting it as a weighted version of the case $q=2$, where the weights themselves become optimization variables.

\section {Source Localization in 2D: SLCP}\label{sec:slcp}

\begin{figure}[tb]
  \centering
  \includegraphics[scale=0.7]{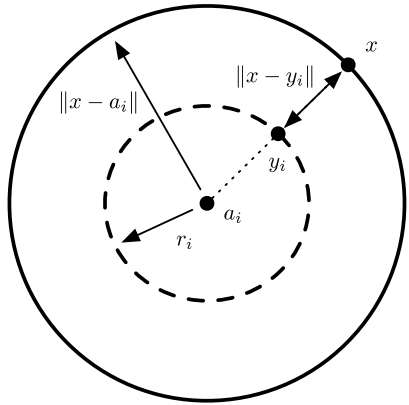}
  \caption{Geometrical interpretation of terms in the source localization cost function \eqref{eq:costfunction1_2} for $p = 1$, $q = 2$.}
  \label{fig:circle}
\end{figure}

For $p = 1$, $q = 2$ we view each term in \eqref{eq:costfunction1_2} as the squared distance between two circles centered on $a_{i}$, one with radius $\|x-a_{i}\|$, and the other with radius $r_{i}$ (see Figure \ref{fig:circle}).
This term can be replaced by the squared norm of the difference between the position vector $x$ and its closest point on the circle $\{y \in \mathbb{R}^{2}: \|y-a_{i}\| = r_{i} \}$, which we denote by $y_{i}$. Problem \eqref{eq:costfunction1_2} can then be equivalently expressed as (a formal proof of equivalence is provided in \cite{Ekim2011})
\begin{equation}\label{sourcecostfunction1}
\begin{array}{cl}
\mbox{minimize}&\sum_{i = 1}^{m} \|x-y_{i}\|^{2}\\
x,y_i&\\
\mbox{subject to}&\|y_{i}-a_{i} \|=r_{i} ~~i = 1, \ldots, m.
\end{array}
\end{equation}
If we fix $y_{i}$, the solution of \eqref{sourcecostfunction1} with respect to $x$ is an unconstrained optimization problem whose solution is readily obtained as the center of mass of the constellation $x = \frac{1}{m}\sum_{i=1}^{m}y_{i}$. Moreover, in 2D the constraints of \eqref{sourcecostfunction1} can be compactly described in the complex plane, yielding 
\begin{equation}\label{sourcecostfunction2}
\begin{array}{cl}
\mbox{minimize}& \|\frac{1}{m} {\bf 1}_m {\bf 1}_m^T\mathbf{y}-\mathbf{y}\|^2\\
\mathbf{y},\bm{\theta}&\\
\mbox{subject to}& \mathbf{y = a + R\bm{\theta}},
\end{array}
\end{equation}
where
$\mathbf{a}=\begin{bmatrix} a_{1} & \ldots & a_{m} \end{bmatrix}^{T} \in \mathbb{C}^{m}$ holds the anchor coordinates, expressed as complex numbers,
$\mathbf{R}= \text{diag}(r_{1}, \ldots, r_{m})
 \in \mathbb{R}^{m \times m}$, and
$\bm{\theta} = \begin{bmatrix} e^{j\phi_1} & \ldots & e^{j\phi_m} \end{bmatrix}^{T} \in \mathbb{C}^{m}$. The complex representation makes it simple to impose unit magnitude constraints on the elements of $\bm{\theta}$, and later relax them to obtain an SDR. Expanding the objective function and deleting constant terms yields the quadratic constrained problem
\begin{equation}\label{sourcecostfunction7}
\begin{array}{cl}
\mbox{minimize}&2\,\text{Re}({\bf c}^H\bm{\theta})-\frac{1}{m}\bm{\theta}^H{\bf rr}^T\bm{\theta}\\
\bm{\theta}\\
\mbox{subject to}&|\theta_i|=1,\\
\end{array}
\end{equation}
where $\mathbf{r} = \mathbf{R}\mathbf{1}_{m}$ and $\mathbf{c} = \mathbf{R}(\mathbf{I}_{m}-\frac{1}{m} \mathbf{1}_{m} \mathbf{1}_{m}^{T}) \mathbf{a}$.

To proceed we now wish to replace $\text{Re}(\mathbf{c}^{H}\bm{\theta})$ in \eqref{sourcecostfunction7} with $-|\mathbf{c}^{H}\bm{\theta}|$, which is readily written as a function of a quadratic form in $\bm{\theta}$ and then relaxed in the same way as the second term in the objective function. To this end, first note that if $\bm{\theta}$ is replaced with $\bm{\theta}e^{j\gamma}$ neither the second term in the objective function of \eqref{sourcecostfunction7} nor the constraints change for any angle $\gamma$. By proper choice of $\gamma$ the complex number $\mathbf{c}^{H}\bm{\theta}$ may be rotated to the (negative) real axis for any feasible $\bm{\theta}$, such that $\text{Re}(\mathbf{c}^{H}\bm{\theta}e^{j\gamma}) = -|\mathbf{c}^{H}\bm{\theta}|$, thus reducing the value of the objective function relative to other values of $\gamma$. This implies that any optimal solution of \eqref{sourcecostfunction7} will satisfy $\text{Re}(\mathbf{c}^{H}\bm{\theta}) = -|\mathbf{c}^{H}\bm{\theta}|$, which justifies replacing $\text{Re}(\cdot)$ with $-|\cdot|$ in the cost function. It should be kept in mind, however, that once a solution $\bm{\theta}$ to the modified optimization problem is obtained it should be rotated to obtain the actual vector of phases $\bm{\theta}e^{j\gamma}$ such that $\text{Re}(\mathbf{c}^{H}\bm{\theta}e^{j\gamma}) = -|\mathbf{c}^{H}\bm{\theta}|$.

Now the modified problem is equivalently written as 
\begin{equation}\label{sourcecostfunction9}
\begin{array}{cl}
\mbox{maximize}&2\sqrt{\text{tr}(\mathbf{c}\mathbf{c}^{H}\bm{\theta}\bm{\theta}^{H})} + \frac{1}{m}\text{tr}(\mathbf{r}\mathbf{r}^{T}\bm{\theta}\bm{\theta}^{H})\\
\bm{\theta}\\
\mbox{subject to}&|\theta_{i}|=1,\\
\end{array}
\end{equation}
and following standard manipulations we introduce the new variable $\bm{\Phi} = \bm{\theta\theta}^H$ and an associated (nonconvex) constraint $\text{rank}(\bm{\Phi}) = 1$. Finally, a SDR formulation of SLCP is obtained by introducing the hypograph variable $t$ such that $0 \leq t \leq 2\sqrt{\text{tr}(\mathbf{c}\mathbf{c}^{H}\bm{\Phi})}$ and dropping the rank constraint 
\begin{equation}\label{eq:costfunction1_2_11}
\begin{array}{cl}
\mbox{maximize}&t+\frac{1}{m}\text{tr}(\mathbf{rr}^T \bm{\Phi})\\
\bm{\Phi},t\\
\mbox{subject to}&\bm{\Phi} \succeq 0, \; \; \phi_{ii}=1, \; \; 4\mathbf{c}^H\bm{\Phi}\mathbf{c}\geq t^2.
\end{array}
\end{equation}
Remark that the solution of \eqref{eq:costfunction1_2_11} is a positive semidefinite matrix, which should have a clearly dominant eigenvalue in problem instances where the SDR is an accurate approximation to the initial problem \eqref{sourcecostfunction1}. In such cases $\bm{\Phi} \approx \lambda_{1}\mathbf{u}_{1}\mathbf{u}_{1}^{H}$, where $\lambda_{1}$ is the highest eigenvalue of $\bm{\Phi}$ and $\mathbf{u}_{1}$ the corresponding eigenvector, and the vector of complex phases is estimated as $\bm{\theta} = \sqrt{\lambda_{1}}\mathbf{u}_{1}$ \cite{Golub1996}. An alternative approach for computing $\bm{\theta}$ is examined in Section \ref{sec:rank1}. Table \ref{tab:slcp} summarizes the SLCP algorithm.
\begin{table}[tb]
 \caption{Summary of the SLCP algorithm}
 \label{tab:slcp}
 \begin{enumerate}
 \item{Given the anchor positions and range measurements, solve the SDR \eqref{eq:costfunction1_2_11}}
 \item{Compute a rank-1 approximation of the SDR solution as $\bm{\Phi} \approx \bm{\theta}\bm{\theta}^{H}$}
\item{Compute a rotation angle $\gamma$ such that $\text{Re}(\mathbf{c}^{H}\bm{\theta}e^{j\gamma}) = -|\mathbf{c}^{H}\bm{\theta}|$ in \eqref{sourcecostfunction7}}
\item{Obtain the vector of circle projections $\mathbf{y} = \mathbf{a} + \mathbf{R}\bm{\theta}e^{j\gamma}$}
\item{Estimate the source position as the centroid $x = \frac{1}{m}\mathbf{1}_{m}^{T}\mathbf{y}$}
 \end{enumerate}
\end{table}


\subsection{Tightness and Geometry of the Constraint Set in SLCP}\label{sec:convexhull}

The source localization problem prior to relaxation \eqref{sourcecostfunction9} can be written as
\begin{equation}\label{sourcecostfunction10}
\begin{array}{cl}
\mbox{maximize}&2\sqrt{u} + \frac{1}{m}v\\
u, v\\
\mbox{subject to}& (u,v) \in {\cal S},
\end{array}
\end{equation}
where
\begin{equation}
  \label{eq:slcpset}
  {\cal S} = \left\{ \left( |\mathbf{c}^{H}\bm{\theta}|^{2}, \; |\mathbf{r}^{T}\bm{\theta}|^{2}\right) : \bm{\theta} \in \mathbb{C}^{m}, \: |\theta_{i}| = 1 \right\}.
\end{equation}
The objective function in \eqref{sourcecostfunction10} is concave with respect to $u$ and $v$, and the optimization problem would be convex if the set ${\cal S}$, over which this function should be maximized, were convex. Then, the SDR used in SLCP \eqref{eq:costfunction1_2_11} would always find a rank-1 solution $\bm{\Phi}$, from which the vector of phases $\bm{\theta}$ would readily follow by factorization. In practice it was found that, even for a moderate number of anchors, the set ${\cal S}$ is likely to have the required shape along part of its border, as discussed below, so that the SDR solution has indeed rank-1. We now examine some of the properties of ${\cal S}$ and the optimal solution.

\begin{figure}[tb]
\centering
\subfloat[three anchors]{\includegraphics[width=7.0cm]{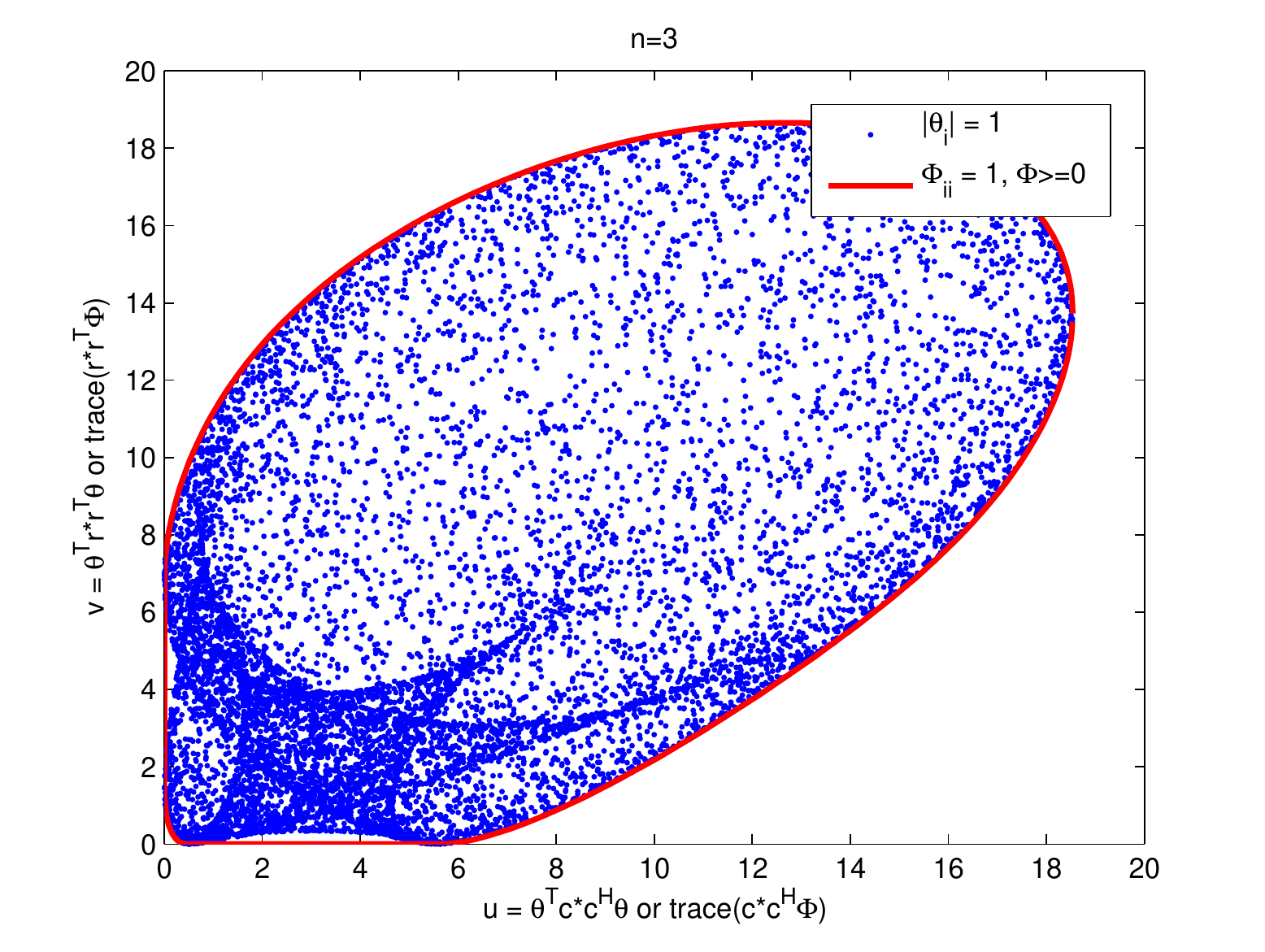}}\hfill
\subfloat[four anchors]{\includegraphics[width=7.0cm]{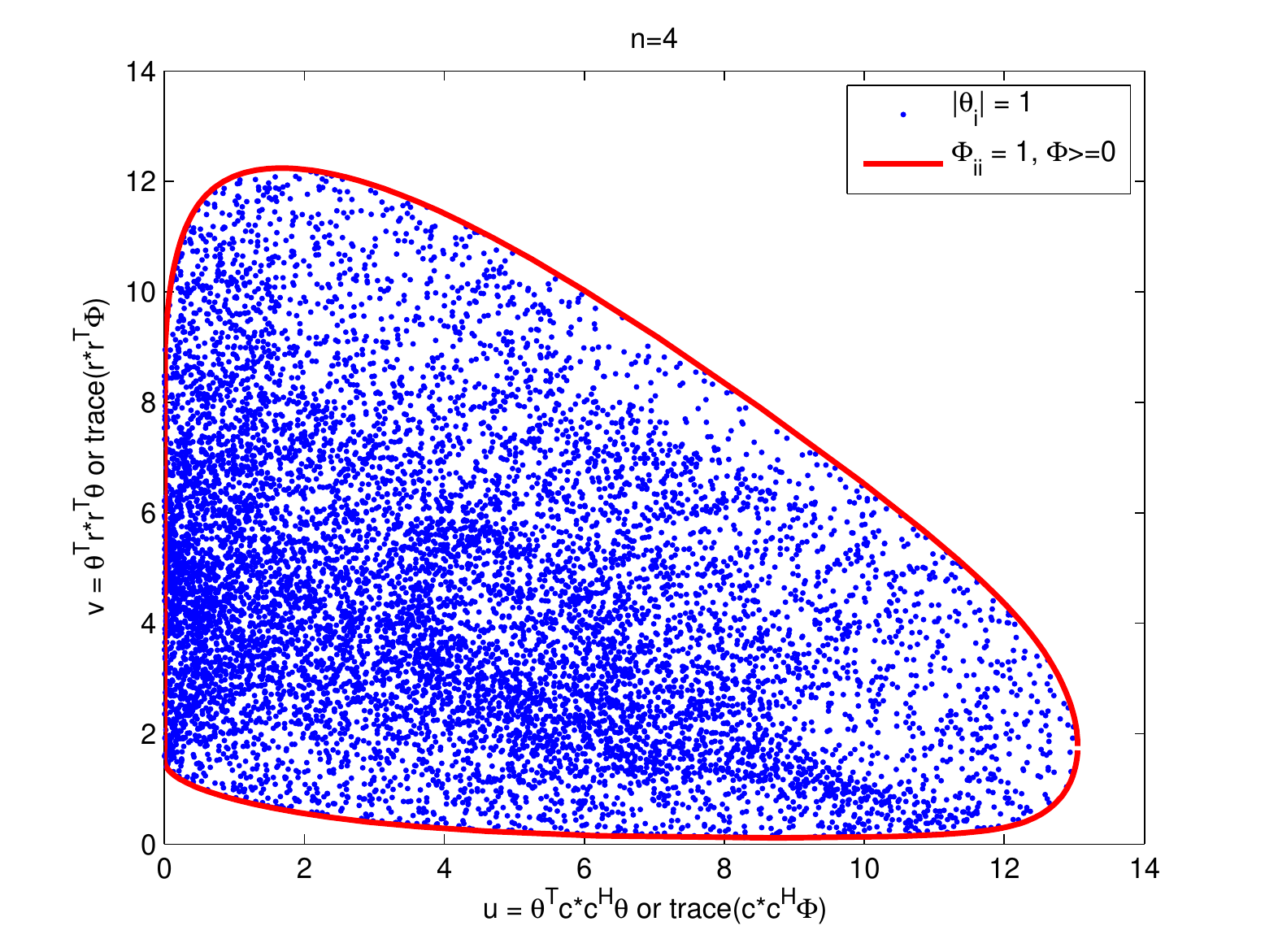}}\\
\subfloat[five anchors]{\includegraphics[width=7.0cm]{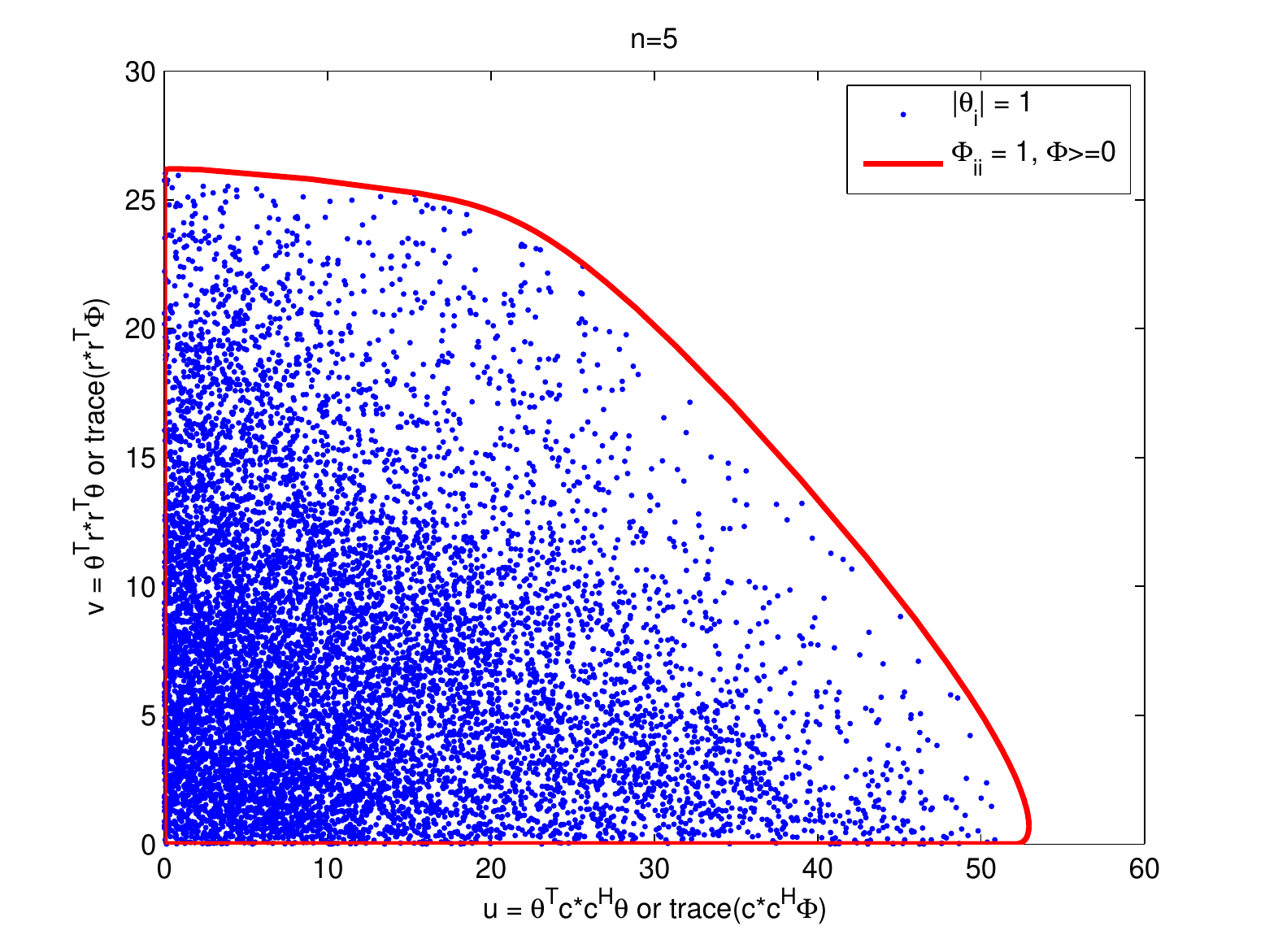}}\hfill
\subfloat[six anchors]{\includegraphics[width=7.0cm]{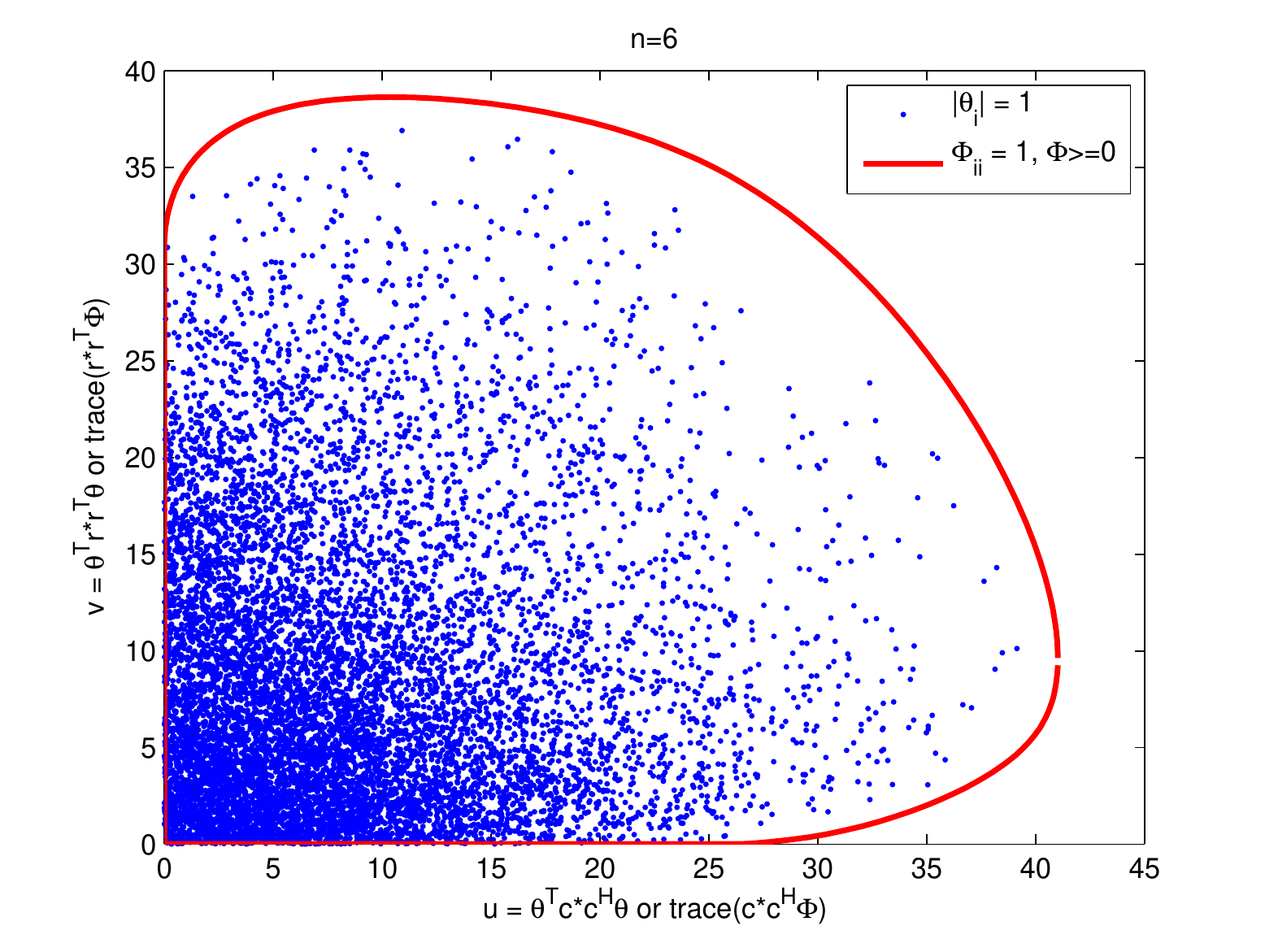}}
\caption{The constraint set ${\cal S}$ for randomly generated $\bm{\theta}$ satisfying $|\theta_i|=1$ and different numbers of randomly placed anchors. For each set the hypothesized convex hull, computed by relaxation of ${\cal S}$, is also depicted.}
\label{fig:convexhull}
\end{figure}

Given the separable form of the cost function \eqref{sourcecostfunction10} it is clear that, for fixed $v$, it can be maximized by choosing $u$ as large as possible within ${\cal S}$, and vice-versa. This implies the following property for the optimal points of \eqref{sourcecostfunction10}:
\begin{prop}
The optimal points of \eqref{sourcecostfunction10} lie on the ``upper right'' boundary of set ${\cal S}$, i.e., optimal points of \eqref{sourcecostfunction10} are maximal elements of ${\cal S}$ with respect to the standard cone $\mathbb{R}^{2}_{+}$ \cite{Boyd2004}.
\end{prop}
Regarding the convexity properties of ${\cal S}$, recall that the cost function of \eqref{sourcecostfunction10} was designed to be invariant to rotations of $\bm{\theta}$ so that, without loss of generality, the first element may be taken as unity. For $m = 2$ anchors and $\theta_{1} = 1$, $\theta_{2} = e^{j\phi}$ we then have
\begin{align}
  u & = |c_{1}^{*}+c_{2}^{*}e^{j\phi}|^{2} = |c_{1}|^{2}+|c_{2}|^{2}+2|c_{1}||c_{2}|\cos(\phi+\alpha) \\
  v & = |r_{1}+r_{2}e^{j\phi}|^{2} = r_{1}^{2}+r_{2}^{2}+2r_{1}r_{2}\cos\phi,
\end{align}
where $\alpha = \angle c_{1} - \angle c_{2}$. Set ${\cal S}$ is an ellipse centered on $(|c_{1}|^{2}+|c_{2}|^{2}, \; r_{1}^{2}+r_{2}^{2})$, therefore clearly nonconvex. Given the definitions of $\mathbf{c}$ and $\mathbf{r}$ in \eqref{sourcecostfunction7}, for $m > 2$ anchors it is always possible to zero out elements $3, \ldots, m$ in these vectors if $r_{3} = \ldots = r_{m} = 0$ in the diagonal of $\mathbf{R}$, thus reverting to the case $m = 2$. In summary:
\begin{prop}
  Depending on the specific range measurements, set ${\cal S}$ may be nonconvex for any number of anchors.
\end{prop}
In spite of the lack of convexity guarantees for ${\cal S}$, our simulation results suggest that for $m \geq 3$ anchors and typical range measurements this set usually does have a convex-like shape. Even when ${\cal S}$ is not convex all that is required for our SDR to provide a rank-1 solution is ``local convexity'' along the ``upper right'' boundary of ${\cal S}$ where the optimal point of \eqref{sourcecostfunction10} is known to be located. More formally, we require that the intersection of ${\cal S}$ with any supporting hyperplane defined by a normal direction with nonnegative components be a compact subset (a single point or a line segment) \cite{Boyd2004}. Figure~\ref {fig:convexhull} depicts some examples of ${\cal S}$ for different numbers of anchors and randomly generated $\mathbf{c}$, $\mathbf{r}$.
Our practical test for local (non)convexity of ${\cal S}$ consists of tracing multiple supporting hyperplanes with nonnegative normal elements, and assessing whether any of them intersect ${\cal S}$ at two well-separated points. We build supporting hyperplanes not on ${\cal S}$ directly, which is a hard problem, but on the related relaxed convex set
\begin{equation}
  \label{eq:slcpcoset}
  {\cal T} = \left\{ \left( \text{tr}(\mathbf{c}\mathbf{c}^{H}\bm{\Phi}), \; \text{tr}(\mathbf{r}\mathbf{r}^{T}\bm{\Phi}) \right) : \bm{\Phi} \in \mathbb{C}^{m\times m}, \: \bm{\Phi} \succeq 0, \: \phi_{ii} = 1 \right\}.
 \end{equation}
Specifically, for a supporting hyperplane with normal $(\cos\beta, \: \sin\beta)$, $0 \leq \beta \leq \frac{\pi}{2}$, we determine an intersection point with ${\cal T}$ by solving the convex optimization problem
\begin{equation}
  \label{eq:supporthyperplane}
\begin{array}{cl}
\mbox{maximize}&\langle (\cos\beta, \, \sin\beta) , \left( \text{tr}(\mathbf{c}\mathbf{c}^{H}\bm{\Phi}), \, \text{tr}(\mathbf{r}\mathbf{r}^{T}\bm{\Phi}) \right) \rangle\\
\bm{\Phi}\\
\mbox{subject to}& \bm{\Phi} \succeq 0, \; \phi_{ii} = 1,
\end{array}
\end{equation}
and setting the intersection point as $u = \text{tr}(\mathbf{c}\mathbf{c}^{H}\bm{\Phi})$, $v = \text{tr}(\mathbf{r}\mathbf{r}^{T}\bm{\Phi})$. This procedure is justified by the following result, proved in Appendix \ref{app:convexhull}.
\begin{thm}\label{thm:convexhull}
For $m \leq 3$ anchors the sets ${\cal S}$ and ${\cal T}$ have the same set of supporting hyperplanes with nonnegative normal elements. Equivalently, in the relevant portion of its boundary ${\cal T}$ coincides with the convex hull of ${\cal S}$.
\end{thm}
Although we only prove this result up to $m = 3$, the empirical evidence suggests that it is also valid for higher $m$, at least up to some maximum order (see Figure \ref{fig:convexhull}). We leave this as a conjecture and apply the procedure for $m > 3$ as well, noting, however, that the case $m = 3$ has major practical significance as the minimum number of anchors that are necessary to recover a general 2D source position based on range measurements. We also conjecture that ${\cal T}$ is actually the convex hull of ${\cal S}$, so \eqref{eq:supporthyperplane} may be used to trace the full boundary of $\text{co}({\cal S})$, and not just the portion where the supporting hyperplanes have nonnegative normal elements. This assumption is not required for our analysis, but was used for generating the set boundaries shown in Figure \ref{fig:convexhull}.


\subsection{Factorization of the SDR Solution}\label{sec:rank1}

The solution of the relaxed SLCP optimization problem (\ref{eq:costfunction1_2_11}) is a positive semidefinite matrix, $\bm{\Phi}$, from which the vector of complex exponentials $\bm{\theta}$ is calculated by rank-1 factorization. The latter is needed to form the vector of circle projections $\mathbf{y} = \mathbf{a} + \mathbf{R}\bm{\theta}$ (see \eqref{sourcecostfunction2}) and, ultimately, the source position vector as the centroid $x = \frac{1}{m}\mathbf{1}_{m}^{T}\mathbf{y}$. The rank-1 factorization method advocated at the end of Section \ref{sec:slcp} is truncation of the eigenvalue decomposition of $\bm{\Phi}$ at the highest eigenvalue. Here we examine a more exact search-based alternative for the practically relevant case of $m = 3$ anchors, which will also be useful to assess the accuracy of the factorization based on eigenvalue truncation.

For a given positive semidefinite matrix $\bm{\Phi} \in \mathbf{C}^{m \times m}$ we wish to find vector $\bm{\theta} \in \mathbf{C}^{m}$ satisfying
\begin{equation}\label{rank10}
\begin{array}{cl}
\mbox{minimize}&\|\bm{\Phi}-\bm{\theta}\bm{\theta}^{H}\|_{F}^{2}\\
\bm{\theta}\\
\mbox{subject to}&|\theta_{i}|=1.\\
\end{array}
\end{equation}
The objective function in \eqref{rank10} is expanded as
\begin{multline}
  \|\bm{\Phi}-\bm{\theta}\bm{\theta}^{H}\|_{F}^{2} = \text{tr}\bigl( (\bm{\Phi}-\bm{\theta}\bm{\theta}^{H})^{H}(\bm{\Phi}-\bm{\theta}\bm{\theta}^{H}) \bigr)
  = \|\bm{\Phi}\|_{F}^{2} + \underbrace{\|\bm{\theta}\|^{4}}_{m^{2}} - \text{tr}(\bm{\Phi}^{H}\bm{\theta}\bm{\theta}^{H})- \text{tr}(\bm{\theta}\bm{\theta}^{H}\bm{\Phi}).
\end{multline}
Ignoring constant terms the optimization problem is equivalently reformulated as
\begin{equation}\label{rank12}
\begin{array}{cl}
\mbox{maximize}&\bm{\theta}^{H}\bm{\Phi}\bm{\theta}\\
\bm{\theta}\\
\mbox{subject to}&|\theta_{i}|=1.\\
\end{array}
\end{equation}
The cost function of \eqref{rank12} is insensitive to a global rotation of all elements of $\bm{\theta}$ by a common factor, hence for $m = 3$ anchors $\bm{\theta}$ can be written as $\bm{\theta} =
\begin{bmatrix} 1 & e^{j\alpha} & e^{j(\alpha+\delta)} \end{bmatrix}^{T}$ and \eqref{rank12} becomes
\begin{equation}\label{rank13}
\begin{array}{cl}
\mbox{maximize}&\text{Re}( \phi_{12}e^{j\alpha}+\phi_{23}e^{j\delta}+\phi_{13}e^{j(\alpha+\delta)} ).\\
\alpha,\delta
\end{array}
\end{equation}
For fixed $\alpha$ the maximum is attained for $\delta = -\angle (\phi_{23}+\phi_{13}e^{j\alpha})$, yielding for \eqref{rank13}
\begin{equation}\label{rank15}
\begin{array}{cl}
\mbox{maximize}&\text{Re}(\phi_{12}e^{j\alpha})+|\phi_{23}+ \phi_{13}e^{j\alpha}|.\\
\alpha&\\
\end{array}
\end{equation}
The solution to \eqref{rank15} is found by searching for the maximum value over the interval $[0, \, 2\pi)$.

Referring to the definitions of the 2D sets ${\cal S}$ in \eqref{eq:slcpset} and ${\cal T}$ in \eqref{eq:slcpcoset}, similar criteria to the above were considered for finding $\bm{\theta}$ such that the induced point in ${\cal S}$ is closest in Euclidean norm to the one induced by $\bm{\Phi}$ in ${\cal T}$. However, the many-to-one nature of the mapping of $\bm{\theta}$ onto points in ${\cal S}$ makes this formulation intrinsically ambiguous.

\section{Source Localization in Higher Dimensions: SLNN}\label{sec:slnn}

To extend the approach used in SLCP to $n > 2$ dimensions, we write the circle/sphere equations in \eqref{sourcecostfunction1} using an equivalent parametric form with real coordinates
\begin{equation}\label{eq:costfunction1_2_nn}
\begin{array}{cl}
\mbox{minimize}&\sum_{i=1}^{m}\|x-y_i\|^{2}\\
x,y_i,u_{i}&\\
\mbox{subject to}&y_i = a_i+r_iu_i, \quad \|u_{i}\| = 1,\\
\end{array}
\end{equation}
where $x, y_i, a_i$ and $u_i$ are now vectors in $\mathbb{R}^{n}$, rather than complex scalars used in SLCP. In \eqref{eq:costfunction1_2_nn} $u_{i} \in \mathbb{R}^{n}$ is a unit-norm vector that plays the same role as the complex phase shift $e^{j\phi_{i}}$ in SLCP. Equivalently,
\begin{equation}\label{eq:costfunction1_2_nn_frobenius}
\begin{array}{cl}
\mbox{minimize}&\|\mathbf{1}_{m}x^{T}-\mathbf{Y}\|_{F}^{2}\\
x,y_i,u_{i}&\\
\mbox{subject to}&\underbrace{\begin{bmatrix} y_{1}^{T}\\ \vdots \\ y_{m}^{T} \end{bmatrix}}_{\mathbf{Y}} = \underbrace{\begin{bmatrix} a_{1}^{T}\\ \vdots \\ a_{m}^{T} \end{bmatrix}}_{\mathbf{A}} + \mathbf{R}\underbrace{\begin{bmatrix} u_{1}^{T}\\ \vdots \\ u_{m}^{T} \end{bmatrix}}_{\mathbf{U}}, \quad \|u_{i}\| = 1,\\
\end{array}
\end{equation}
where $\mathbf{R}= \text{diag}(r_{1}, \ldots, r_{m})$ as in \eqref{sourcecostfunction2}. For fixed $y_{i}$, $u_{i}$ \eqref{eq:costfunction1_2_nn_frobenius} describes $n$ uncoupled least-squares problems whose variables are the components of the source location vector $x$. The optimal solutions may be jointly written compactly as
\begin{equation}
  x^{T} = (\mathbf{1}_{m}^{T}\mathbf{1}_{m})^{-1}\mathbf{1}_{m}^{T}\mathbf{Y} = \frac{1}{m}\mathbf{1}_{m}^{T}\mathbf{Y}.
\end{equation}
Replacing this back in \eqref{eq:costfunction1_2_nn_frobenius} to eliminate variable $x$ the objective function becomes $\| \bm{\Pi}\mathbf{Y} \|_{F}^{2} = \text{tr}(\mathbf{Y}^{T}\bm{\Pi}\mathbf{Y})$, where $\bm{\Pi} = \mathbf{I}_{m} - \frac{1}{m}\mathbf{1}_{m}\mathbf{1}_{m}^{T}$ is a projection matrix (hence idempotent). Similarly to \eqref{sourcecostfunction2}--\eqref{sourcecostfunction7} we can now eliminate variable $\mathbf{Y}$ and the first set of equality constraints, expanding its definition in the objective function and ignoring constant terms to obtain
\begin{equation}\label{eq:costfunction1_2_nn_3}
\begin{array}{cl}
\mbox{minimize}& 2\,\text{tr}(\mathbf{C}^{T}\mathbf{U})-\frac{1}{m}\text{tr}(\mathbf{ U}^{T}\mathbf{rr}^{T}\mathbf{U})\\
\mathbf{U}&\\
\mbox{subject to}&\|u_{i}\|=1,\\
\end{array}
\end{equation}
where $\mathbf{C} = \mathbf{R}\bm{\Pi}\mathbf{A}$ and, as in \eqref{sourcecostfunction7}, $\mathbf{r} = \mathbf{R}\mathbf{1}_{m}$.


\paragraph{Nuclear Norm Approximation}


As in the complex formulation we wish to rewrite the first term in the objective function of \eqref{eq:costfunction1_2_nn_3} in a form that is more amenable to SDR. In the optimization problem we thus replace $\mathbf{U}$ with the product $\mathbf{U}\mathbf{V}$, where $\mathbf{V}$ is an $n \times n$ orthogonal matrix such that $\mathbf{V}^{T}\mathbf{V} = \mathbf{V}\mathbf{V}^{T} = \mathbf{I}_{n}$, yielding
\begin{equation}\label{eq:costfunction1_2_nn_4}
\begin{array}{cl}
\mbox{minimize}& 2\,\text{tr}(\mathbf{C}^{T}\mathbf{U}\mathbf{V})-\frac{1}{m}\text{tr}(\mathbf{V}^{T}\mathbf{ U}^{T}\mathbf{rr}^{T}\mathbf{U}\mathbf{V})\\
\mathbf{U}, \mathbf{V}&\\
\mbox{subject to}&\|u_{i}\|=1, \quad \mathbf{V}^{T}\mathbf{V} = \mathbf{I}_{n}.\\
\end{array}
\end{equation}
Note that, due to the orthogonality of $\mathbf{V}$, each line of $\mathbf{U}\mathbf{V}$ still has unit norm, so for any feasible $\mathbf{U}$ in \eqref{eq:costfunction1_2_nn_3} $\mathbf{U}\mathbf{V}$ is also feasible. Regarding \eqref{eq:costfunction1_2_nn_4}, $\mathbf{V}$ may be interpreted as an inner optimization variable that, for each candidate $\mathbf{U}$, minimizes the value of the objective function. Noting that the second term in the objective function \eqref{eq:costfunction1_2_nn_4} does not depend on $\mathbf{V}$, as $\text{tr}(\mathbf{V}^{T}\mathbf{ U}^{T}\mathbf{rr}^{T}\mathbf{U}\mathbf{V}) = \text{tr}(\mathbf{rr}^{T}\mathbf{U}\mathbf{V}\mathbf{V}^{T}\mathbf{ U}^{T}) = \text{tr}(\mathbf{rr}^{T}\mathbf{U}\mathbf{ U}^{T})$, the inner optimization problem simply becomes
\begin{equation}\label{eq:costfunction1_2_nn_4_inner}
\begin{array}{cl}
\mbox{minimize}& \text{tr}(\mathbf{C}^{T}\mathbf{U}\mathbf{V}) = \langle \mathbf{V}, \mathbf{U}^{T}\mathbf{C} \rangle\\
\mathbf{V}&\\
\mbox{subject to}&\mathbf{V}^{T}\mathbf{V} = \mathbf{I}_{n}.\\
\end{array}
\end{equation}
This involves the minimization of a linear function on the set of orthogonal matrices, which resembles the known problem of minimizing a linear function of a vector $\mathbf{v}$, say, $\langle \mathbf{v},\mathbf{a} \rangle$, on the unit sphere $\|\mathbf{v}\|^{2} = \mathbf{v}^{T}\mathbf{v} = 1$. Invoking the KKT conditions \cite{Boyd2004} the latter problem is readily seen to yield the optimal cost $-\|\mathbf{a}\|$, attained at the point on the sphere along vector $-\mathbf{a}$. One would therefore expect the solution of \eqref{eq:costfunction1_2_nn_4_inner} to be $-\|\mathbf{C}^{T}\mathbf{U}\|$, involving some matrix norm of $\mathbf{C}^{T}\mathbf{U}$. In Appendix \ref{app:slnn} it is shown that this is indeed the case, and that the appropriate norm to consider is the \emph{nuclear norm}, defined for matrix $\mathbf{X}$ as $\|\mathbf{X}\|_{N} = \text{tr}\bigl( (\mathbf{X}^{H}\mathbf{X})^{\frac{1}{2}}\bigr)$, and equaling the sum of its singular values \cite{Petersen2008}. The optimization problem \eqref{eq:costfunction1_2_nn_4} is therefore equivalently rewritten as
\begin{equation}\label{eq:costfunction1_2_nn_7}
\begin{array}{cl}
\mbox{minimize}& -2\|\mathbf{C}^{T}\mathbf{U}\|_{N} - \frac{1}{m}\text{tr}(\mathbf{ rr}^{T}\mathbf{U}\mathbf{U}^{T})\\
\mathbf{U}&\\
\mbox{subject to}&\|u_{i}\| = 1,
\end{array}
\end{equation}
or
\begin{equation}\label{eq:costfunction1_2_nn_8}
\begin{array}{cl}
\mbox{maximize}& 2\,\text{tr}\bigl( (\mathbf{C}^{T}\mathbf{U}\mathbf{U}^{T}\mathbf{C})^{\frac{1}{2}} \bigr) + \frac{1}{m}\text{tr}(\mathbf{rr}^{T}\mathbf{U}\mathbf{U}^{T})\\
\mathbf{U}&\\
\mbox{subject to}&\|u_{i}\| = 1.
\end{array}
\end{equation}
We now introduce the variable $\mathbf{W} = \mathbf{U}\mathbf{U}^{T}$ and ignore the associated nonconvex constraint $\text{rank}(\mathbf{W}) = n$ to obtain the SDR
\begin{equation}\label{eq:costfunction1_2_nn_9}
\begin{array}{cl}
\mbox{maximize}& 2\,\text{tr}\bigl( (\mathbf{C}^{T}\mathbf{W}\mathbf{C})^{\frac{1}{2}} \bigr) + \frac{1}{m}\text{tr}(\mathbf{rr}^T\mathbf{W})\\
\mathbf{W}&\\
\mbox{subject to}&\mathbf{W}\succeq 0, \quad w_{ii} = 1.\\
\end{array}
\end{equation}
The objective function of \eqref{eq:costfunction1_2_nn_9} is the sum of a concave\footnote{The first term is the composition of the linear map $\mathbf{X} = \mathbf{C}^{T}\mathbf{W}\mathbf{C}$ with $\text{tr}(\mathbf{X}^{\frac{1}{2}})$, which is known to be concave in $\mathbf{X}$ \cite{Boyd2004}.} function of $\mathbf{W}$ with a linear term, and is therefore concave. The constraint set of \eqref{eq:costfunction1_2_nn_9} is convex, thus establishing that this is indeed a convex optimization problem. We express it in standard SDP form as
\begin{equation}\label{eq:costfunction1_2_nn_11}
\begin{array}{cl}
\mbox{maximize}& 2\,\text{tr}(\mathbf{Z}) + \frac{1}{m}\text{tr}(\mathbf{rr}^{T}\mathbf{W})\\
\mathbf{W,Z}&\\
\mbox{subject to}&\mathbf{W}\succeq 0, \quad w_{ii}=1, \quad
\begin{bmatrix} \mathbf{C}^{T}\mathbf{W}\mathbf{C} & \mathbf{Z}\\ \mathbf{Z} & \mathbf{I}_{n} \end{bmatrix} \succeq 0, \quad \mathbf{Z} \succeq 0.
\end{array}
\end{equation}
The equivalence between \eqref{eq:costfunction1_2_nn_9} and \eqref{eq:costfunction1_2_nn_11} is proved in Appendix \ref{app:slnn}.

Similarly to the complex 2D formulation, the solution of our SDR is a $m \times m$ matrix $\mathbf{W}$ that should have approximately rank $n$ when the relaxation is tight. The matrix $\mathbf{U}$ of unit-norm vectors is obtained by SVD factorization of $\mathbf{W}$ \cite{Golub1996} and, after accounting for the inner rotation of $\mathbf{U}$, it is used to build the $y_{i}$ and, ultimately, the source position vector $x$. Table \ref{tab:slnn} summarizes the SLNN algorithm.
\begin{table}[tb]
 \caption{Summary of the SLNN algorithm}
 \label{tab:slnn}
 \begin{enumerate}
 \item{Given the anchor positions and range measurements, solve the SDR \eqref{eq:costfunction1_2_nn_11}}
 \item{Compute a rank-$n$ approximation of the SDR solution as $\mathbf{W} \approx \mathbf{U}\mathbf{U}^{T}$}
\item{Solve the inner optimization problem \eqref{eq:costfunction1_2_nn_4_inner} to get the rotation matrix $\mathbf{V}$}
\item{Obtain the matrix of sphere projections as $\mathbf{Y} = \mathbf{A} + \mathbf{R}\mathbf{U}\mathbf{V}$}
\item{Estimate the source position as the centroid of the rows of $\mathbf{Y}$, $x = \frac{1}{m}\mathbf{Y}^{T}\mathbf{1}_{m}$}
 \end{enumerate}
\end{table}


\subsection{Localization under Laplacian Noise: SL-$\ell_1$}\label{sec:SLl1}


When disturbances are Laplacian and i.i.d., thus heavier tailed than Gaussian, maximizing the likelihood amounts to solving \eqref{eq:costfunction1_2} for $p=q=1$,
\begin{equation}\label{sourcecostfuncl1}
\begin{array}{cl}\mbox{minimize}& \sum_{i=1}^{m} |\|x-a_i\|-r_{i}|.\\
x
\end{array}
\end{equation}
The presence of $\abs{\cdot}$ in each summation term of \eqref{sourcecostfuncl1}, rather than $(\cdot)^{2}$, de-emphasizes the contributions of measurements $r_{i}$ corrupted by large noise values. The optimal point of \eqref{sourcecostfuncl1} is thus less biased by these outlier measurements than the cost function \eqref{eq:costfunction1_2} for the Gaussian case $p=1$, $q=2$. However, a major difficulty in solving \eqref{sourcecostfuncl1} is the fact that the cost function is not differentiable, making it less amenable to the types of analytic manipulations that we use to develop SDR. The strategy that we adopt to circumvent this difficulty parallels the one used in \cite{Ekim2011} for 2D sources, and as a key ingredient involves squaring the cost function of \eqref{sourcecostfuncl1} (which does not affect the location of extremal points), and then  rewriting it as
\begin{equation}\label{sourcecostfuncl1.4a}
\begin{array}{cl}
\mbox{minimize}&\sum_{i=1}^{m} \frac{(\|x-a_{i}\|-r_{i})^{2}}{\lambda_{i}}\\
x,\bm{\lambda}\\
\text{subject to} & \lambda_{i}>0, \quad \mathbf{1}_{m}^{T}\bm{\lambda}=1.
\end{array}
\end{equation}
The cost function is thus reduced to a weighted version of the more tractable Gaussian log-likelihood, where the real weighting coefficients $\lambda_{i}$ become optimization variables themselves. See \cite{Ekim2011} for a proof of this result (also \cite{Nemirovski2001}). Now, the manipulations used earlier in Section \ref{sec:slnn} for the development of SLNN can be replicated here to reformulate the problem as
\begin{equation}\label{sourcecostfuncl1.4}
\begin{array}{cl}
\mbox{minimize}&\sum_{i=1}^{m} \frac{\|x-y_{i}\|^{2}}{\lambda_{i}}\\
x,y_{i},u_{i},\bm{\lambda}\\
\text{subject to} & y_{i} = a_{i}+r_{i}u_{i}, \quad \|u_{i}\| = 1, \quad \lambda_{i}>0, \quad \mathbf{1}_{m}^{T}\bm{\lambda}=1.
\end{array}
\end{equation}
For given $y_{i}$, $u_{i}$, and $\bm{\lambda}$, \eqref{sourcecostfuncl1.4} has a least-squares cost function whose unconstrained optimal solution with respect to $x$ is readily found in closed form from the first-order stationarity condition
\begin{align}
  \sum_{i=1}^{m}\frac{x-y_{i}}{\lambda_{i}} & = 0, & x^{*} & = \frac{\sum_{i=1}^{m}\frac{y_{i}}{\lambda_{i}} }{\sum_{i=1}^{m}\frac{1}{\lambda_{i}}}.
\end{align}
Substituting the optimal $x$ in \eqref{sourcecostfuncl1.4}, and using matrix notation, the cost function becomes $\text{tr}(\mathbf{Y}^{T}\bm{\Xi}\mathbf{Y})$, where $\bm{\Xi}$ is the modified projector
\begin{equation}\label{eq:proj_lambda}
  \begin{split}
\bm{\Xi} &= \left[\begin{array}{ccc}
\frac{1}{\lambda_1}& & 0\\ 
& \ddots &\\
0 & & \frac{1}{\lambda_m}\end{array}\right]
-\frac{1}{\sum_{i=1}^{m}\frac{1}{\lambda_i}}
\left [
\begin{array}{c}
\frac{1}{\lambda_1}\\
\vdots\\
\frac{1}{\lambda_m} \end{array}\right]
\left [
\begin{array}{ccc}
\frac{1}{\lambda_1}&\ldots&\frac{1}{\lambda_m} \end{array}\right] = \bm{\Lambda}^{-1} - \bm{\Lambda}^{-1}\mathbf{1} (\mathbf{1}^{T}\bm{\Lambda}^{-1}\mathbf{1})^{-1}\mathbf{1}^{T}\bm{\Lambda}^{-1}
 \end{split}
\end{equation}
with $\bm{\Lambda} = \text{diag}(\lambda_1,\ldots,\lambda_{m})$.

\subsubsection{Alternating directions (SL-$\ell_1$ AD)}

One possibility for iteratively solving \eqref{sourcecostfuncl1.4} is to use block coordinate descent, alternating between minimizing the expression with respect to $\{x, \, y_{i}, \, u_{i}\}$ for fixed $\bm{\lambda}$ and vice-versa. For fixed $\bm{\lambda}$ the problem is 
\begin{equation}\label{eq:costfunction1_1_3D_2}
\begin{array}{cl}
\mbox{minimize}& \text{tr}\bigl( (\mathbf{A}+\mathbf{R}\mathbf{U})^{T}\bm{ \Xi}(\mathbf{A}+\mathbf{R}\mathbf{U}) \bigr)\\
\mathbf{U}\\
\mbox{subject to}&\|u_{i}\|=1,\\
\end{array}
\end{equation}
which differs from the SLNN formulation only in the projector matrix $\bm{\Xi}$. It can therefore be similarly manipulated into a relaxed form that parallels \eqref{eq:costfunction1_2_nn_11}
\begin{equation}\label{eq:costfunction1_1_3D_4}
\begin{array}{cl}
\mbox{maximize}& 2\,\text{tr}(\mathbf{Z}) + \frac{1}{\kappa}\text{tr}(\mathbf{rr}^{T}\mathbf{W})\\
\mathbf{W,Z}&\\
\mbox{subject to}&\mathbf{W}\succeq 0, \quad w_{ii}=1, \quad \begin{bmatrix} \mathbf{C}^{T}\mathbf{W}\mathbf{C} & \mathbf{Z}\\ \mathbf{Z} & \mathbf{I}_{n} \end{bmatrix} \succeq 0, \quad \mathbf{Z} \succeq 0,
\end{array}
\end{equation}
with $\mathbf{C} = \mathbf{R}\bm{\Xi}\mathbf{A}$, $\mathbf{r} = \mathbf{R} \begin{bmatrix} \frac{1}{\lambda_1} & \ldots & \frac{1}{\lambda_m} \end{bmatrix}^{T}$, and $\kappa = \sum_{i=1}^{m}\frac{1}{\lambda_{i}}$. For the converse block coordinate descent step with fixed $\{x, \, y_{i}, \, u_{i}\}$ the problem is
\begin{equation}\label{eq:costfunction1_1_3D_coordstep_lambda}
\begin{array}{cl}
\mbox{minimize}&\sum_{i=1}^{m} \frac{K_{i}^{2}}{\lambda_{i}}\\
\bm{\lambda}\\
\text{subject to} & \lambda_{i}>0, \quad \mathbf{1}_{m}^{T}\bm{\lambda}=1,
\end{array}
\end{equation}
where $K_{i} \stackrel{\Delta}{=} |\|x-a_{i}\|-r_{i}| = \|x-y_{i}\|$ are constant in this subproblem. The solution, readily obtained from the first-order KKT conditions, is given by \cite{Ekim2011}
\begin{equation}\label{eq:update_lambda}
  \lambda_{i}^{*} = \frac{K_{i}}{\sum_{i=1}^{m}K_{i}} = \frac{ |\|x-a_{i}\|-r_{i}|}{\sum_{i=1}^{m} |\|x-a_{i}\|-r_{i}|},
\end{equation}
yielding the desired $\ell_{1}$-type cost function $\sum_{i=1}^{m} \frac{K_{i}^{2}}{\lambda_{i}} = (\sum_{i=1}^{m} K_{i})^{2}$. We sequentially perform the iterations \eqref{eq:costfunction1_1_3D_4} and \eqref{eq:update_lambda}, starting with $\bm{\lambda} = \frac{1}{m}\mathbf{1}_{m}$, until $\|x^{k+1}-x^{k}\|$ is within some prescribed tolerance $\varepsilon$. In our simulated scenarios on the order of 3--10 iterations are needed for $\varepsilon = 10^{-2}$. This method is denoted by SL-$\ell_1$ AD.

\subsubsection{Non-iterative formulation (SL-$\ell_1$ MD)}

For a non-iterative solution of \eqref{sourcecostfuncl1.4} we start from the equivalent formulation \eqref{eq:costfunction1_1_3D_2}, with $\lambda_{1}, \ldots, \lambda_{m}$ included as optimization variables through the weighting matrix $\bm{\Xi}$, and introduce an epigrapth variable $t_{i}$ for each term contributing to $\text{tr}(\cdot)$ in the cost function
\begin{equation}\label{eq:costfunction1_1_3D_epigraph}
\begin{array}{cl}
\mbox{minimize}& \mathbf{t}\mathbf{1}_{n} \\
\mathbf{U}, \bm{\lambda}, \mathbf{t} \\
\mbox{subject to}& \mathbf{e}_{i}^{T}(\mathbf{A}+\mathbf{R}\mathbf{U})^{T}\bm{ \Xi}(\mathbf{A}+\mathbf{R}\mathbf{U})\mathbf{e}_{i} \leq t_{i} \\
& \|u_{i}\|=1, \quad \lambda_{i}>0, \quad \mathbf{1}_{m}^{T}\bm{\lambda} = 1,
\end{array}
\end{equation}
where $\mathbf{t} = \begin{bmatrix} t_{1} & \ldots & t_{n} \end{bmatrix}$ and $\mathbf{e}_{i}$ is the standard coordinate vector with 1 in the $i$-th position and zeros elsewhere. As in \cite{Ekim2011} we invoke the matrix inversion lemma to express \eqref{eq:proj_lambda} as the limiting case of (positive semidefinite) $\bm{\Xi} = \lim_{\sigma \to \infty}(\bm{\Lambda} + \sigma \mathbf{1}_{m}\mathbf{1}_{m}^{T})^{-1}$, which is more amenable to analytic manipulations in optimization problems. In practice we take $\sigma$ as a sufficiently large constant. Using Schur complements the inequality constraint in \eqref{eq:costfunction1_1_3D_epigraph} may be successively written as
\begin{align}
  \begin{bmatrix}
    t_i & \mathbf{e}_i^T(\mathbf{A+RU})^T\\
    (\mathbf{A+RU})\mathbf{e}_i & \bm{\Xi}^{-1}
  \end{bmatrix} & \succeq 0 \label{eq:schur_l1_1} \\
t_{i}(\bm{\Lambda} + \sigma \mathbf{1}_{m}\mathbf{1}_{m}^{T}) - (\mathbf{A+RU})\mathbf{e}_{i} \mathbf{e}_{i}^{T}(\mathbf{A+RU})^{T} & \succeq 0. \label{eq:schur_l1_2}
\end{align}
The last inequality is bilinear in $t_{i}$ and $\lambda_{1}, \ldots, \lambda_{m}$, and we linearize it by replacing the optimization variable $\bm{\lambda}$ with a new $\bm{\beta}_{i} = t_{i}\bm{\lambda}$. Now, the $\bm{\beta}_{i}$ can be assembled into a matrix
\begin{equation}\label{eq:assemble_beta}
  \bm{\beta} = \begin{bmatrix} \bm{\beta}_{1} & \ldots & \bm{\beta}_{n} \end{bmatrix} = \bm{\lambda}\mathbf{t},
\end{equation}
which, as shown above, should have rank 1 and satisfy $\beta_{ij} > 0$, $\mathbf{1}_{m}^{T}\bm{\beta} = \mathbf{t}$. However, the rank-1 constraint for $\bm{\beta}$ cannot be directly imposed in convex formulations, and we resort to a common technique to indirectly induce low rank in optimal solutions by adding to the cost function the (scaled) nuclear norm $\|\bm{\beta}\|_{N}$.

Regarding the second term on the left-hand side of \eqref{eq:schur_l1_2}, we first note that
\begin{equation}\label{eq:factor_colsU}
  (\mathbf{A+RU})\mathbf{e}_{i} = \begin{bmatrix} \mathbf{A}\mathbf{e}_{i} & \mathbf{R} \end{bmatrix} \begin{bmatrix}1 \\ \mathbf{U}\mathbf{e}_{i} \end{bmatrix} = \begin{bmatrix} \bm{\alpha}_{i} & \mathbf{R} \end{bmatrix} \begin{bmatrix} 1 \\ \bm{\upsilon}_{i} \end{bmatrix},
\end{equation}
where $\bm{\alpha}_{i}$ and $\bm{\upsilon}_{i}$ denote the $i$-th columns of matrices $\mathbf{A}$ and $\mathbf{U}$, respectively. Now, consider the following variable, obtained from the stacked rotation vectors that make up $\mathbf{U}$,
\begin{equation}\label{eq:relaxU}
  \begin{split}
    \mathbf{W} & = \begin{bmatrix} 1 \\ \text{vec}(\mathbf{U}^T) \end{bmatrix} \begin{bmatrix} 1 & \text{vec}(\mathbf{U}^T)^{T} \end{bmatrix}  =
    \begin{bmatrix}
      1 & \mathbf{u}_{1}^{T} & \ldots & \mathbf{u}_{m}^{T} \\
      \mathbf{u}_{1} & \underbrace{\mathbf{u}_{1}\mathbf{u}_{1}^{T}}_{\mathbf{W}_{11}} \\
      \vdots & & \ddots \\
      \mathbf{u}_{m} & & & \underbrace{\mathbf{u}_{m}\mathbf{u}_{m}^{T}}_{\mathbf{W}_{mm}}
    \end{bmatrix}.
\end{split}
\end{equation}
Further, let $I_{i}$ denote the set of row indices that extracts the elements of $[1 \; \bm{\upsilon}_{i}^{T}]^{T}$ in \eqref{eq:factor_colsU} from the first column of $\mathbf{W}$. Then, the dyad below is readily obtained by selecting the submatrix formed from the $I_{i}$ rows and $I_{i}$ columns of $\mathbf{W}$
\begin{equation}
  \mathbf{W}_{I_{i}I_{i}} = \begin{bmatrix} 1 \\ \bm{\upsilon}_{i} \end{bmatrix} \begin{bmatrix} 1 & \bm{\upsilon}_{i}^{T} \end{bmatrix},
\end{equation}
and this carries over to \eqref{eq:schur_l1_2} through \eqref{eq:factor_colsU}, which can therefore be written in terms of submatrix $\mathbf{W}_{I_{i}I_{i}}$. The positive semidefinite matrix $\mathbf{W}$ will replace $\mathbf{U}$ as an  optimization variable, retaining the constraints along the diagonal blocks in \eqref{eq:relaxU}, namely, $\text{tr}(\mathbf{W}_{ii}) = 1$. Finally, we obtain the full convex relaxation of \eqref{eq:costfunction1_1_3D_epigraph} by combining all the above elements and dropping the rank-1 constraint for $\mathbf{W}$ that is implied by \eqref{eq:relaxU}
\begin{equation}\label{varT3}
\begin{array}{cl}
\mbox{minimize}& \mathbf{t}\mathbf{1}_{n} + \mu\|\bm{\beta}\|_{N}\\
\mathbf{W}, \bm{\beta}, \mathbf{t} \\
\mbox{subject to}& \text{diag}(\bm{\beta}_{i}) + t_{i} \sigma \mathbf{1}_{m}\mathbf{1}_{m}^{T} \succeq \begin{bmatrix} \bm{\alpha}_{i} & \mathbf{R} \end{bmatrix} \mathbf{W}_{I_{i}I_{i}} \begin{bmatrix} \bm{\alpha}_{i}^{T} \\ \mathbf{R} \end{bmatrix} \\
& \mathbf{W} \succeq 0, \quad w_{11} = 1, \quad \text{tr}(\mathbf{W}_{ii}) = 1, \quad  \beta_{ij}>0, \quad \mathbf{1}_{m}^{T}\bm{\beta} = \mathbf{t}.
\end{array}
\end{equation}
This reference formulation for SL-$\ell_{1}$ in multiple dimensions is denoted by SL-$\ell_{1}$ MD.

\subsubsection{Simplified non-iterative formulation (SL-$\ell_1$ SD)}

Our simulation results suggest that in most scenarios the accuracy of the solution obtained from \eqref{varT3} is nearly identical to that of a simplified formulation where a single epigraph variable, $t$, is used. Referring to \eqref{eq:costfunction1_1_3D_epigraph}, we now minimize $\text{tr}(t\mathbf{I}_{n})$ or, equivalently, $t$, and replace the first constraint for all $i = 1, \ldots, n$ with the single matrix inequality $(\mathbf{A}+\mathbf{R}\mathbf{U})^{T}\bm{ \Xi}(\mathbf{A}+\mathbf{R}\mathbf{U}) \preceq t\mathbf{I}_{n}$. Applying Schur complements as in \eqref{eq:schur_l1_1}--\eqref{eq:schur_l1_2} yields
\begin{equation}
t(\bm{\Lambda} + \sigma \mathbf{1}_{m}\mathbf{1}_{m}^{T}) - \begin{bmatrix} \mathbf{A} & \mathbf{R}\end{bmatrix} \begin{bmatrix} \mathbf{I}_{n} \\ \mathbf{U}\end{bmatrix} \begin{bmatrix} \mathbf{I}_{n} & \mathbf{U}^{T}\end{bmatrix} \begin{bmatrix} \mathbf{A}^{T} \\ \mathbf{R}\end{bmatrix} \succeq 0,
\end{equation}
and again we replace variable $\bm{\lambda}$ with $\bm{\beta} = t\bm{\lambda}$ such that $\beta_{i} > 0$, $\mathbf{1}_{m}^{T}\bm{\beta} = t$. Now, however, there is no need to assemble a matrix as in \eqref{eq:assemble_beta} and to include its nuclear norm as a penalization term in the cost function. Finally, to obtain a convex relaxation we replace $\mathbf{U}$ with the new variable
\begin{equation}\label{eq:relaxUsing}
  \mathbf{W} = \begin{bmatrix} \mathbf{I}_{n} \\ \mathbf{U}\end{bmatrix} \begin{bmatrix} \mathbf{I}_{n} & \mathbf{U}^{T}\end{bmatrix} = \begin{bmatrix} \underbrace{\mathbf{I}_{n}}_{\mathbf{W}_{11}} & \mathbf{U}^{T} \\ \mathbf{U} & \mathbf{U}\mathbf{U}^{T} \end{bmatrix},
\end{equation}
and drop the rank-$n$ constraint on $\mathbf{W}$ that follows from \eqref{eq:relaxUsing}. The simplified SDP formulation for SL-$\ell_{1}$ in multiple dimensions, denoted by SL-$\ell_{1}$ SD, is given by
\begin{equation}\label{eq:singT}
\begin{array}{cl}
\mbox{minimize}& t \\
\mathbf{W}, \bm{\beta}, t \\
\mbox{subject to}& \text{diag}(\bm{\beta}) + t\sigma \mathbf{1}_{m}\mathbf{1}_{m}^{T} \succeq \begin{bmatrix} \mathbf{A} & \mathbf{R}\end{bmatrix} \mathbf{W} \begin{bmatrix} \mathbf{A}^{T} \\ \mathbf{R}\end{bmatrix} \\
& \mathbf{W} \succeq 0, \quad \mathbf{W}_{11} = \mathbf{I}_{n}, \quad w_{ii} = 1, \quad \beta_{i}>0, \quad \mathbf{1}_{m}^{T}\bm{\beta} = t.
\end{array}
\end{equation}
Note that the optimization variables $\mathbf{W}$ and $\bm{\beta}$ in \eqref{varT3} have size $(mn+1) \times (mn+1)$ and $m \times n$, respectively, whereas the corresponding sizes in \eqref{eq:singT} are only $(m+n) \times (m+n)$ and $m \times 1$. For ambient dimension $n = 2$ or 3 and for $m \approx 5$ anchors used in our simulations problem \eqref{eq:singT} has considerably fewer variables than \eqref{varT3}, and the gap increases as $m$ and $n$ grow.

Given the configuration for variable $\mathbf{W}$ in both non-iterative formulations of SL-$\ell_{1}$ \eqref{eq:relaxU}, \eqref{eq:relaxUsing}, the required elements of the rotation vectors that make up $\mathbf{U}$ can be obtained from the rightmost (block) column of $\mathbf{W}$ or by factorizing submatrices along the block diagonal. The former approach is usually more accurate \cite{Ye2004}.


\section{Numerical Results}\label{sec:simu}


In this section the tightness and the accuracy of our source localization algorithms are tested in 2D and 3D scenarios, and under various noise assumptions. Results are benchmarked against another relaxation-based method proposed in \cite{Cheung2004}, denoted below as SDR (as in \cite{Stoica2008}), and with the Squared Range LS (SR-LS) approach of \cite{Stoica2008}. While the latter does not resort to relaxation, but rather directly optimizes the source coordinates using an iterative root-finding procedure, we take its performance as representative of the current state of the art in optimization-based source localization.

In each reported simulation we performed $M$ Monte Carlo runs, where in each run the source and anchor locations were randomly generated from a uniform distribution over a square or cube whose sides are $[-10,10]$. The observed ranges, corrupted by i.i.d.\ noise, were generated as described in Section~\ref{sec:formulation} under appropriate noise probability densities. The tables list Root Mean-Square Errors (RMSE), computed as $\sqrt{\frac{1}{M}\sum_{i=1}^{M} \|x_{i} - \hat{x}_{i}\|^{2}}$, where $x_{i}$ and $\hat{x}_{i}$ denote the actual and estimated source positions in the $i$-th Monte Carlo run, respectively. For ease of reference, the best result among all algorithms tested for any given setup will often be shown in boldface.

\vspace{1ex}\noindent \textbf{\emph{Convexity and tightness of SLCP:}} In this example we characterize the accuracy of the convex relaxation used in SLCP and compare its performance to that of the SDR algorithm of \cite{Cheung2004}. Range measurements to a variable number of randomly placed anchors were generated as indicated above over $M = 1000$ Monte Carlo runs, and corrupted by white Gaussian noise.

First, we estimate how often the constraint set $\cal S$ \eqref{eq:slcpset}, which appears in our formulation of the source localization problem prior to relaxation \eqref{sourcecostfunction10}, is convex along its ``upper right'' boundary where the optimal solution lies. As discussed in Section \ref{sec:convexhull}, when this property holds the relaxed solution $\bm{\Phi}$ obtained by SLCP \eqref{eq:costfunction1_2_11} will have rank 1 and can be factorized to yield the optimal point for the non-relaxed problem \eqref{sourcecostfunction10} on the boundary of $\cal S$. We empirically assess convexity of $\cal S$ by tracing the boundary of the (partially hypothesized) convex hull $\cal T$ \eqref{eq:slcpcoset} and searching for line segments that delimit regions where the boundaries of $\cal S$ and $\cal T$ depart due to local concavity of $\cal S$.  Specifically, we solve the support hyperplane problem \eqref{eq:supporthyperplane} for a grid of angles $0 \leq \beta \leq \frac{\pi}{2}$ and detect the presence of a line segment when the distance between the intersection points $\bigl( u(\beta),v(\beta) \bigr)$ for two consecutive angles $\beta$ exceeds a threshold. For a noise standard deviation $\sigma_{\text{gaussian}} = 10^{-2}$, $\cal S$ passed the convexity test in 80\% of runs for three anchors. The percentage increased to 84\% for five anchors, in line with our reasoning in Section \ref{sec:convexhull} that $\cal S$ is more likely to be convex as the number of anchors increases.

Next, we compare the RMSEs of SDR and SLCP. As in \cite{Beck2008} we provide results for all Monte Carlo runs (denoted by SDR, SLCP) and also for so-called \emph{tight runs} (denoted by SDRt, SLCPt) where the solution for the relaxed localization problem is close to having rank 1, as desired for subsequent factorization to obtain the actual source coordinates. We consider a solution matrix to be tight when the ratio between its first and second eigenvalues is at least $10^2$. Table \ref{tab:tightrun} lists the RMSEs and the number of tight runs ($N_{\text{SDR}}$, $N_{\text{SLCP}}$) over 1000 trials for five anchors and Gaussian noise standard deviations of $1$, $10^{-1}$, $10^{-2}$, and $10^{-3}$. SLCP is clearly superior over the full set of trials, but the gap to SDR closes in the subset of tight runs, indicating that the advantage is mostly due to a much higher probability of its solution having near rank 1. Even for the highest noise power, where the number of tight runs in both algorithms is comparable, the ratio of first to second eiganvalues is usually higher in SLCP, leading to lower RMSE.
\begin{table}[htb]
\caption{Source localization accuracy for relaxation-based methods (RMSEs listed for total and tight runs).}
\centering
\begin{tabular}{|c|c|c|c|c|c|c|c|c|} \hline
$\sigma_{\text{gaussian}}$ & $N_{\text{SDR}}$ & $N_{\text{SLCP}}$ & SDR \cite{Cheung2004} & SDRt \cite{Cheung2004} & SLCP & SLCPt\\ \hline
$10^{-3}$ &490 & 921 & 0.0045 & 0.0014 & 0.0020 & 0.0015\\
$10^{-2}$ &444 & 815 & 0.0162 & 0.0107 & 0.0112 & 0.0108\\ 
$10^{-1}$ &478 & 527 & 0.1503 & 0.0960 & 0.1207 & 0.0959\\ 
1      &538 & 526 & 1.6070 & 1.1885 & 1.2169 & 1.1885\\ \hline \hline
\end{tabular}
\label{tab:tightrun}
\end{table}

Under the same simulation setup as above, but using only three anchors, we test the alternative search-based method described in Section \ref{sec:rank1} to obtain the vector of rotation factors $\bm{\theta}$ from the relaxed solution matrix of SLCP, $\bm{\Phi}$. Improvements in total RMSE are under 1\% for all noise variances using $2\times10^5$ grid points on the interval $[0, \, 2\pi)$ to evaluate \eqref{rank15}. Foremost, this suggests that rank-1 factorization by SVD, which we adopt as our technique of choice to efficiently extract rotation factors, yields results that are indeed very close to the best possible strategy for finding $\bm{\theta}$.


\vspace{1ex}\noindent \textbf{\emph{Localization in 2D and 3D under Gaussian noise:}} In the remaining simulations our algorithms are benchmarked against SR-LS \cite{Stoica2008}, whose global performance exceeds that of SDR \cite{Cheung2004} because it directly optimizes over source locations and does not experience degradations related to tightness of the solutions in the same way that SDR does.

Measurements for the 2D case were generated as in the convexity/tightness assessment, using five anchors and $M = 200$ Monte Carlo runs. Table \ref{tab:rmse_nn} lists the RMSEs for SR-LS, SLCP \eqref{eq:costfunction1_2_11} and its multidimensional counterpart SLNN \eqref{eq:costfunction1_2_nn_11}, and also the algorithms for Laplacian noise, namely, the complex formulation of SL-$\ell_{1}$ in \cite{Ekim2010}, and its multidimensional counterpart SL-$\ell_{1}$ AD \eqref{eq:costfunction1_1_3D_4}, \eqref{eq:update_lambda}.
\begin{table}[htb]
\caption{\label{tab:rmse_nn}2D source localization under Gaussian noise (5 anchors, 200 Monte Carlo runs). RMSEs are given for complex (SLCP, SL-$\ell_{1}$) and real (SLNN, SL-$\ell_{1}$ AD) formulations.}
\centering
\begin{tabular}{|c|c|c|c|c|c|c|c|}\hline \hline
$\sigma_{\text{gaussian}}$ & SR-LS \cite{Stoica2008} & SLCP & SL-$\ell_1$ \cite{Ekim2010} & SLNN & SL-$\ell_1$ AD\\ \hline
$10^{-3}$ &0.0032 & 0.0023 & $\bm{0.0014}$ & 0.0023 & 0.0057\\
$10^{-2}$ &0.0138 & $\bm{0.0109}$ & 0.0136 & 0.0113 & 0.0133\\ 
$10^{-1}$ &0.1406 & $\bm{0.1037}$ & 0.1118& 0.1097 & 0.1249\\ 
1       &1.4947 & $\bm{1.3249}$ & 1.4536 & 1.3580 & 1.4593\\ \hline \hline
\end{tabular}
\end{table}
All algorithms outperform SR-LS, which squares measurements ($p = 2$, $q = 2$ in \eqref{eq:costfunction1_2}) and thus becomes more sensitive to the presence of (Gaussian) noise in range measurements. The fact that SLCP/SLNN achieve the best results is not surprising, as these algorithms actually maximize a Gaussian likelihood function. Interestingly, SLCP attains slighly lower errors than SLNN, even though the same cost function and similar steps are used in the derivation of both algorithms. However, the impact of relaxing the rank constraint in the optimization variable to obtain a semidefinite program is not necessarily the same, which could explain the observed differences in performance. Similar comments apply to SL-$\ell_{1}$ and SL-$\ell_{1}$ AD, although the algorithmic differences between the complex and real formulations are larger than for SLCP/SLNN.

Results for 3D source localization, for which the complex formulations cannot be used, are given in Table \ref{tab:rmse_nn_3}.
\begin{table}[htb]
\caption{\label{tab:rmse_nn_3}3D source localization under Gaussian noise (5 anchors, 200 Monte Carlo runs).}
\centering
\begin{tabular}{|c|c|c|c|c|c|}\hline \hline
$\sigma_{\text{gaussian}}$ & SR-LS \cite{Stoica2008} & SLNN & \multicolumn{3}{|c|}{SL-$\ell_1$} \\
\cline{4-6}
& & & AD & MD & SD\\ \hline
$10^{-3}$ & 0.0040 & $\bm{0.0036}$ & 0.0038 & 0.0038 & 0.0040 \\ 
$10^{-2}$ & 0.0295 & $\bm{0.0274}$ & 0.0285 & 0.0290 & 0.0292 \\ 
$10^{-1}$ & 0.2612 & $\bm{0.2290}$ & 0.2376 & 0.2401 & 0.2390 \\
1               & 3.3279 & $\bm{2.7431}$ & 2.9492 & 2.8748 & 2.8801 \\ \hline \hline
\end{tabular}
\end{table}
Again, both SLNN ans SL-$\ell_{1}$ outperform SR-LS, the former having lower RMSE as its cost function is matched to the noise statistics. Source localization with the same number of anchors is a less constrained problem in 3D than it is in 2D, resulting in higher RMSEs. Regarding the three variants of SL-$\ell_{1}$ (see Section \ref{sec:SLl1}), there is no well-defined trend on their relative performance. Note also how the RMSE of the non-adaptive simplified formulation SL-$\ell_{1}$ SD \eqref{eq:singT} is quite close to that of the general formulation SL-$\ell_{1}$ MD \eqref{varT3}, even outperforming it, on average, for one of the noise powers.

\vspace{1ex}\noindent \textbf{\emph{Localization in 2D and 3D in the presence of outliers:}}
The same setup for Gaussian noise is adopted here, except that ranges are contaminated either by Laplacian noise, or by what we designate as \emph{selective Gaussian noise}. Range measurements for the latter are created as $r_{i} = \|x-a_i\|+w_{i}+|\epsilon|$, where $w_{i}$ is a Gaussian noise term with $\sigma_{\text{gaussian}} = 0.04$ that is present in all observations and $\epsilon$ is also a Gaussian disturbance, but with higher standard deviation $\sigma_{\text{outlier}}\in [0.3,1.5]$, that contaminates only one measured range (i.e., $\epsilon = 0$ for all other observations). This statistical model is less tractable than the Laplacian noise model, but we include it in some of our simulations as it more realistically reflects how outliers occur in real ranging systems.

Tables \ref{tab:mse_outlier} and \ref{tab:mse_outlier_3} list RMSEs for 2D and 3D source localization under both outlier generation models.
\begin{table}[t]
\caption{\label{tab:mse_outlier}2D source localization in the presence of outlier range measurements (5 anchors, 200 Monte Carlo runs).}
\centering
\subfloat[Laplacian noise\label{tab:mse_outlier:laplace}]{
\begin{tabular}{|c|c|c|c|c|c|}\hline \hline
$\sigma_{\text{laplacian}}$ & SR-LS \cite{Stoica2008} & SLCP & SL-$\ell_1$ \cite{Ekim2010} & SLNN & SL-$\ell_1$ AD\\ \hline              
0.2 & 1.1398 &  1.0839  &  $\bm{0.9240}$ & 1.0770  & 1.0687\\
0.4 &1.9031 & 1.8585  &  $\bm{1.4546}$  &   1.8908  &1.7996\\ 
0.8 & 3.1543  &  3.1143   &  $\bm{3.0344}$  &  3.0812  &  3.0798  \\ \hline \hline
\end{tabular}}
\hspace{2em}
\subfloat[Selective Gaussian noise\label{tab:mse_outlier:selectgauss}]{
\begin{tabular}{|c|c|c|c|c|c|c|}\hline \hline
$\sigma_{\text{outlier}}$ & SR-LS \cite{Stoica2008} & SLCP & SL-$\ell_1$ \cite{Ekim2010} & SLNN & SL-$\ell_1$ AD\\ \hline
0.5 &  0.2983 &  0.2448 &  $\bm{0.1849}$ &  0.2556 &  0.2337\\
1.0 &  0.4662 &  0.4561 &  $\bm{0.2508}$ &  0.4516 &  0.3714\\
1.5 &  1.2419 &  1.1640 &  $\bm{1.0542}$&  1.2389&  1.2157\\ \hline \hline
\end{tabular}}
\end{table}
\begin{table}[t]
\caption{\label{tab:mse_outlier_3}3D source localization in the presence of outlier range measurements (5 anchors, 200 Monte Carlo runs).}
\centering
\subfloat[Laplacian noise\label{tab:mse_outlier:laplace_3}]{
\begin{tabular}{|c|c|c|c|c|c|}\hline \hline
$\sigma_{\text{laplacian}}$ & SR-LS \cite{Stoica2008} & SLNN & \multicolumn{3}{|c|}{SL-$\ell_1$} \\
\cline{4-6}
& & & AD & MD & SD\\ \hline
0.25 & 1.3619 & 1.3577 & 1.2113 & 1.2097 & $\bm{1.1776}$  \\
0.5   & 2.6514 & 2.5719 & $\bm{2.3236}$ & 2.4704 & 2.3651 \\ 
0.75 & 3.5968 & 3.5070 & $\bm{3.1265}$ & 3.2173 & 3.1987 \\ \hline \hline
\end{tabular}}
 \hspace{2em}
 \subfloat[Selective Gaussian noise\label{tab:mse_outlier:selectgauss_3}]{
 \begin{tabular}{|c|c|c|c|c|c|}\hline \hline
 $\sigma_{\text{outlier}}$ & SR-LS \cite{Stoica2008} & SLNN & \multicolumn{3}{|c|}{SL-$\ell_1$} \\
 \cline{4-6}
 & & & AD & MD & SD\\ \hline
 0.3 & 0.3930 & 0.3237 & $\bm{0.3033}$ &0.3035 & $\bm{0.3033}$ \\
 0.6 & 0.9756 & 0.9154 & $\bm{0.8786}$ & 0.9084 & 0.9001 \\
 0.9 & 1.5934 & 1.5502 & $\bm{1.2755}$ & 1.2857 & 1.3022 \\ \hline \hline
 \end{tabular}}
\end{table}
When compared with Tables \ref{tab:rmse_nn} and \ref{tab:rmse_nn_3} for the Gaussian case, the most striking difference is that the variants of SL-$\ell_{1}$, designed for Laplacian noise, now outperform both SR-LS and SLCP/SLNN. As in the 2D Gaussian case, the complex formulation SL-$\ell_{1}$ attains lower errors than the real formulation SL-$\ell_{1}$ AD. In 3D (Table \ref{tab:mse_outlier_3}) the real variants of SL-$\ell_{1}$ developed in Section \ref{sec:SLl1} provide significant RMSE reductions, on the order of 10\%, over those of the Gaussian algorithm SLNN. Interestingly, this conclusion still holds for the selective Gaussian case, whose outlier generation model is not matched to the Laplacian assumption underlying SL-$\ell_{1}$. Similarly to the Gaussian case of Table \ref{tab:rmse_nn_3}, the three variants of SL-$\ell_{1}$ exhibit similar performance, here with a slight advantage of the alternating direction algorithm SL-$\ell_{1}$ AD.

\section{Conclusion}\label{sec:conclusion}

We have proposed SLNN as an extention to 3 and higher dimensions of the ML-based source localization approach developed in \cite{Ekim2010} by formulating it as an optimization problem using nuclear norms and SDR. Similarly, we also extended to higher dimensions the 2D SL-$\ell_1$ localization algorithm for Laplacian noise developed in \cite{Ekim2011}. Our simulation results show that the proposed algorithms provide very accurate results compared to other optimization-based localization methods that operate on range measurements, although their performance in 2D is not quite as good as that of the complex formulations developed in \cite{Ekim2010,Ekim2011}. In 3D scenarios with Gaussian noise SLNN delivered solutions that were about 5\% more accurate than those of SL-$\ell_{1}$, whereas in the presence of outlier range measurements the situation was reversed and SL-$\ell_{1}$ proved to be about 5--10\% more accurate under either Laplacian or selective Gaussian models. We developed both iterative (SL-$\ell_{1}$ AD) and non-iterative (SL-$\ell_{1}$ MD/SD) 3D extensions of SL-$\ell_{1}$, which exhibited comparable performance, with a slight advantage of SL-$\ell_{1}$ AD. Complexity considerations (e.g., computational load, maximum admissible problem size) will then play an important role when selecting one of those algorithms for a particular application.

We have carried out an analysis of the geometry of our 2D formulation for ML localization under Gaussian noise (SLCP), and found that the high probability that a certain portion of the (outer) border of its constraint set is convex justifies the observed strong tightness of our relaxation. Our simulation results for random anchor configurations indicate that another well-known SDR relaxation for the same problem has a significantly higher chance of yielding optimal solutions that do not have the necessary properties (unit rank) to accurately recover source positions. Regarding the extraction of spatial coordinates from the positive semidefinite matrix computed by SLCP, we examined a search-based alternative to standard rank-1 factorization using the SVD. This strategy is feasible for the practically important case of range-based localization using three anchors, but was found to yield only minor improvements relative to the SVD-based factorization.

\appendices


\section{Proof of \lemma\ \ref{thm:convexhull}}\label{app:convexhull}


Our proof of \lemma\ \ref{thm:convexhull} relies on a result, interesting in its own right, that characterizes the convex hull of the set of $3\times 3$ rank-1 matrices built from complex vectors with unit-magnitude components.
\begin{thm}\label{thm:convexhulldyads}
Let
\begin{align}
  {\cal A} & = \left\{ \bm{\theta}\bm{\theta}^{H} : \bm{\theta} \in \mathbb{C}^{3}, \: |\theta_{i}| = 1 \right\}, \label{eq:dyadset} \\
  {\cal B} & = \left\{ \mathbf{\Phi} \in \mathbb{C}^{3\times 3} : \bm{\Phi} \succeq 0, \: \phi_{ii} = 1 \right\}.
\end{align}
then ${\cal B} = \text{co}({\cal A})$.
\end{thm}
\begin{IEEEproof}
$\text{co}({\cal A}) \subset {\cal B}$ is straightforward since ${\cal B}$ is convex and ${\cal A} \subset {\cal B}$. For the reverse direction $\text{co}(A) \supset {\cal B}$ our goal is to find, for every $\mathbf{\Phi} \in {\cal B}$, matrices $\bm{\Theta}_{i} \in {\cal A}$ and nonzero scalars $\lambda_{i} \geq 0$, with $\sum_{i} \lambda_{i} = 1$, such that $\bm{\Phi} = \sum_{i} \lambda_{i}\bm{\Theta}_{i}$.

Note that both ${\cal A}$ and ${\cal B}$ are invariant under the (unitary) similarity operation
\begin{equation}
 \mathbf{M} \rightarrow \mathbf{P}\mathbf{M}\mathbf{P}^{H},
\end{equation}
where $\mathbf{P}$ is the product of a permutation and a diagonal unitary matrix. In other words, we can simultaneously permute rows and columns and multiply the $i$-th row and $i$-th column by a unit-magnitude complex number. Thus, we can assume without loss of generality that $\bm{\Phi}$ is of the form
\begin{align}
  \bm{\Phi} & =
  \begin{bmatrix} 1 & a & b\\ a & 1 & z^{*}\\ b & z & 1 \end{bmatrix}, & 0 & \leq a \leq b, \quad z \in \mathbb{C}.
  \label{eq:convhullmtx}
\end{align}
Since $\bm{\Phi} \succeq 0$, we must have $a \leq 1$, $b \leq 1$, $\abs{z} \leq 1$ and
\begin{equation}
  0 \leq |\bm{\Phi}| = 1 - a^{2} - b^{2} -\abs{z}^{2} + 2ab\,\text{Re}\{ z \},
\end{equation}
which, for $z = x + jy$, reads
\begin{equation}
  \label{eq:convexhullcircle}
  (x-ab)^{2} + y^{2} \leq (1-a^{2})(1-b^{2}).
\end{equation}
For fixed $a$, $b$ this inequality describes a circle (with interior) in the $(x, \, y)$ plane, centered on $(ab, \, 0)$. Since any point in the interior of a circle can be written as a convex combination of two points on its boundary, we can assume that we have equality in \eqref{eq:convexhullcircle}. Thus, from now on we assume
\begin{equation}
  z = ab + \sqrt{(1-a^{2})(1-b^{2})}e^{j\varphi}.
  \label{eq:zsdp}
\end{equation}
We now complete the proof by expressing such $\bm{\Phi}$ as a convex combination of two matrices from ${\cal A}$. For given $0 \leq a \leq b \leq 1$ and $\varphi \in [0, \, 2\pi [$ we want to find $\alpha, \, \beta, \, \gamma, \, \delta \in [0, \, 2\pi [$, and $0 \leq \lambda \leq 1$ such that
\begin{equation}
  \begin{split}
    \bm{\Phi} & =
    \begin{bmatrix} 1 & a & b\\ a & 1 & z^{*}\\ b & z & 1 \end{bmatrix} =
    \lambda \begin{bmatrix} 1\\ e^{j\alpha}\\ e^{j\beta} \end{bmatrix}\begin{bmatrix} 1 & e^{-j\alpha} & e^{-j\beta} \end{bmatrix}  \quad + (1-\lambda) \begin{bmatrix} 1\\ e^{j\gamma}\\ e^{j\delta} \end{bmatrix}\begin{bmatrix} 1 & e^{-j\gamma} & e^{-j\delta} \end{bmatrix}.
  \end{split}
  \label{eq:dyaddecomp}
\end{equation}
We thus have
\begin{align}
  a & = \lambda e^{j\alpha} + (1-\lambda) e^{j\gamma}, \quad b = \lambda e^{j\beta} + (1-\lambda) e^{j\delta}, \label{eq:abdyads} \\
z & = \lambda e^{j(\beta-\alpha)} + (1-\lambda) e^{j(\delta-\gamma)}.
\end{align}
From the first two relations we get
\begin{align}
  e^{j\gamma} & = \frac{a - \lambda e^{j\alpha}}{1-\lambda}, & e^{j\delta} & = \frac{b - \lambda e^{j\beta}}{1-\lambda},
\end{align}
and replacing these in the third relation yields, after simple manipulations,
\begin{equation}
  z = ab + \frac{\lambda}{1-\lambda}(e^{-j\alpha}-a)(e^{j\beta}-b).
  \label{eq:zdyads}
\end{equation}
Before proceeding, we state and prove a useful lemma from elementary geometry:
\begin{thm}\label{thm:chords}
  Referring to Figure \ref{fig:chords_lemma}, if $A$ is a point inside a unit circle whose distance to the center is $a$, $RS$ is any line through $A$, and $PQ$ is a diameter through $A$, then
  \begin{equation}
   AR \cdot AS = AP \cdot AQ = (1-a)(1+a) = 1-a^{2}.
  \end{equation}
  \begin{figure}[tb]
    \centering
    \subfloat[]{%
      \label{fig:chords_lemma}%
      \includegraphics[width=0.3\columnwidth]{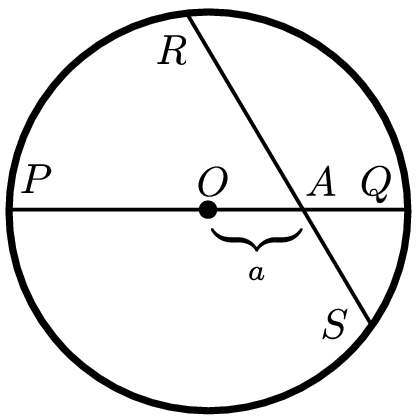}} \hfill
    \subfloat[]{%
      \label{fig:chords_triang}%
      \includegraphics[width=0.3\columnwidth]{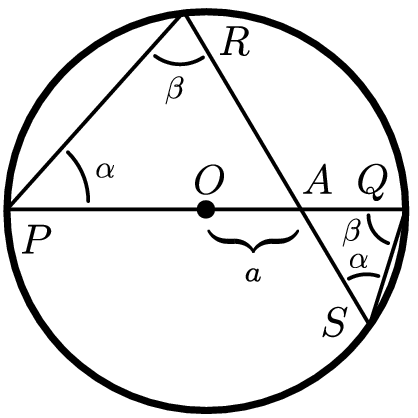}}
    \caption{Illustration of geometrical \lemma\ \ref{thm:chords}.}
    \label{fig:chords}
  \end{figure}
\end{thm}
\begin{IEEEproof}
  Triangles $APR$ and $AQS$, depicted in Figure \ref{fig:chords_triang}, are similar, hence
  \begin{equation}
    \frac{AP}{AS} = \frac{AR}{AQ}.
  \end{equation}
\end{IEEEproof}
We use the lemma above with parameters as depicted in Figure \ref{fig:convcombcircle}.
\begin{figure}[tb]
  \centering
  \includegraphics[width=0.3\columnwidth]{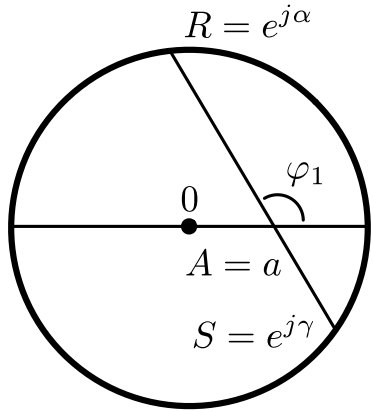}
  \caption{Application of \lemma\ \ref{thm:chords} to a convex combination on the unit circle.}
  \label{fig:convcombcircle}
\end{figure}
From $A = \lambda R + (1-\lambda) S$ we have $\frac{AR}{AS} = \frac{1-\lambda}{\lambda}$, and by \lemma\ \ref{thm:chords} $AR \cdot AS = 1-a^{2}$, hence
\begin{align}
  AR & = \sqrt{\frac{1-\lambda}{\lambda}(1-a^{2})}, & e^{j\alpha} & = a + \sqrt{\frac{1-\lambda}{\lambda}(1-a^{2})}e^{j\varphi_{1}}.
\label{eq:ejalpha}
\end{align}
Similarly, with $A = b$, $R = e^{j\beta}$, $S = e^{j\delta}$, and $\varphi_{2}$ instead of $\varphi_{1}$, we have
\begin{equation}
  e^{j\beta} = b + \sqrt{\frac{1-\lambda}{\lambda}(1-b^{2})}e^{j\varphi_{2}}.
\label{eq:ejbeta}
\end{equation}
Substituting \eqref{eq:ejalpha}, \eqref{eq:ejbeta} back in \eqref{eq:zdyads} yields
\begin{equation}
  z = ab + \sqrt{(1-a^{2})(1-b^{2})}e^{j(\varphi_{2}-\varphi_{1})},
\end{equation}
which has the same form as \eqref{eq:zsdp}, obtained from the positive semidefinite condition for matrix $\bm{\Phi}$ in \eqref{eq:convhullmtx}.

We now argue that letting angle $\alpha$ go from $0$ to $2\pi$ is equivalent to letting $\varphi_{1}$ cover an interval of length $2\pi$ as well (Figure \ref{fig:convcombcircle}). Fixing $\varphi_{1}$, and consequently $\alpha$, the two relations in \eqref{eq:abdyads}, together with an arbitrary requirement that $\text{Im}\{ e^{j\beta} \} \geq 0$, fix the values\footnote{Equivalently, note that fixing $\varphi_{1}$ fully defines the geometrical construction shown in Figure \ref{fig:convcombcircle}, and thus fixes the values of $\gamma$ and $\lambda$. Then, $\lambda$ fully defines the corresponding construction for $A = b$ if, in addition, $\text{Im}\{ e^{j\beta} \} \geq 0$ is specified, and thus fixes the values of $\beta$, $\delta$, and $\varphi_{2}$.} of $\beta$, $\gamma$, $\delta$, $\lambda$, and, in particular, of $\varphi_{2}$. Thus, $\varphi_{2} = f(\varphi_{1})$ is a continuous function of $\varphi_{1}$.

When $\varphi_{1} = 0$, $\varphi_{2}$ has a certain value, say, $\varepsilon_{0} \in [0,\pi]$ (it can be computed, but is not needed in this proof). For $\varphi_{1} = \pi$ it is straightforward to see that $\varphi_{2} = \pi - \varepsilon_{0}$, and for $\varphi_{1} = 2\pi$ it is again $\varepsilon_{0}$. In particular the continuous function $\varphi_{2}-\varphi_{1}$ takes values from $\varepsilon_{0} - 0 = \varepsilon_{0}$ to $\varepsilon_{0} - 2\pi$, i.e., modulo $2\pi$ it takes all values in $[0, \, 2\pi[$. Thus, for any given angle $\varphi$ in \eqref{eq:zsdp}, let $\varphi_{1}$ be such that $f(\varphi_{1})-\varphi_{1} = \varphi$, modulo $2\pi$. Then, the corresponding $\alpha$, $\beta$, $\gamma$, $\delta$, and $\lambda$, as explained above, give the desired decomposition \eqref{eq:dyaddecomp}.
\end{IEEEproof}

We now proceed and prove \lemma\ \ref{thm:convexhull} under the assumption of \lemma\ \ref{thm:convexhulldyads}, thus tacitly assuming $m = 3$. Note that the proof is valid for arbitrary $\mathbf{c}$, $\mathbf{r}$ in \eqref{eq:slcpset} and \eqref{eq:slcpcoset}, i.e., it does not require that the structure for these vectors defined in \eqref{sourcecostfunction7} be taken into account.

\begin{IEEEproof}
  We rewrite sets ${\cal S}$ in \eqref{eq:slcpset} and ${\cal T}$ in \eqref{eq:slcpcoset} using the notation \eqref{eq:dyadset}
\begin{align}
  {\cal S} & = \left\{ \left( \mathbf{c}^{H}\bm{\Theta}\mathbf{c}, \; \mathbf{r}^{T}\bm{\Theta}\mathbf{r} \right) : \bm{\Theta} \in {\cal A} \right\}, \\
  {\cal T} & = \left\{ \left( \mathbf{c}^{H}\bm{\Phi}\mathbf{c}, \; \mathbf{r}^{T}\bm{\Phi}\mathbf{r} \right) : \bm{\Phi} \in \text{co}({\cal A}) \right\}.
\end{align}
Obviously ${\cal S} \subset {\cal T}$. Now let $\alpha \in [0, \, \frac{\pi}{2}]$ and define
\begin{equation}
 (u_{1}, \, v_{1}) = \arg\max_{(u,v) \in {\cal T}} \langle (\cos\alpha, \, \sin\alpha), \, (u, \, v) \rangle.
 \label{eq:supportmaxdpco}
\end{equation}
We wish to show that
\begin{equation}
  \langle (\cos\alpha, \, \sin\alpha), \, (u_{1}, \, v_{1}) \rangle = \max_{(u,v) \in {\cal S}} \langle (\cos\alpha, \, \sin\alpha), \, (u, \, v) \rangle,
  \label{eq:supportmaxdp}
\end{equation}
so that the inner product over ${\cal S}$ attains the same maximum value as over the larger set ${\cal T}$, and the support hyperplanes with normal $(\cos\alpha, \, \sin\alpha)$ thus coincide for the two sets. It is enough to prove that there exists $(u', \, v') \in {\cal S}$ that attains the left-hand side of \eqref{eq:supportmaxdp}.

We may write $\bm{\Phi}_{1} \in \text{co}({\cal A})$ which maximizes \eqref{eq:supportmaxdpco} as
\begin{align}
  \bm{\Phi}_{1} & = \sum_{i} \lambda_{i}\bm{\theta}_{i}\bm{\theta}_{i}^{H}, & \lambda_{i} & \geq 0, \; \sum_{i} \lambda_{i} = 1, \; \abs{\theta_{ik}} = 1,
\end{align}
hence
\begin{equation}
  \begin{split}
  \langle ( & \cos\alpha, \, \sin\alpha), \, (u_{1}, \, v_{1}) \rangle \\
& = (\sqrt{\cos\alpha}\,\mathbf{c})^{H}\bm{\Phi}_{1}(\underbrace{\sqrt{\cos\alpha}\,\mathbf{c}}_{\mathbf{p}}) + (\sqrt{\sin\alpha}\,\mathbf{r})^{T}\bm{\Phi}_{1}(\underbrace{\sqrt{\sin\alpha}\,\mathbf{r}}_{\mathbf{q}}) \\
& = \sum_{i} \lambda_{i} \left( \mathbf{p}^{H}\bm{\theta}_{i}\bm{\theta}_{i}^{H}\mathbf{p} + \mathbf{q}^{H}\bm{\theta}_{i}\bm{\theta}_{i}^{H}\mathbf{q} \right) \\
& = \sum_{i} \lambda_{i} \left( \abs{\mathbf{p}^{H}\bm{\theta}_{i}}^{2} + \abs{\mathbf{q}^{H}\bm{\theta}_{i}}^{2} \right).
 \end{split}
\end{equation}
Let $i_{0}$ be the index where the last summation attains its maximum value. Then
\begin{multline}
 \langle ( \cos\alpha, \, \sin\alpha), \, (u_{1}, \, v_{1}) \rangle \leq \abs{\mathbf{p}^{H}\bm{\theta}_{i_{0}}}^{2} + \abs{\mathbf{q}^{H}\bm{\theta}_{i_{0}}}^{2}
= \langle ( \cos\alpha, \, \sin\alpha), \, (\mathbf{c}^{H}\bm{\theta}_{i_{0}}\bm{\theta}_{i_{0}}^{H}\mathbf{c}, \, \mathbf{r}^{T}\bm{\theta}_{i_{0}}\bm{\theta}_{i_{0}}^{H}\mathbf{r}) \rangle,
\end{multline}
which completes the proof because the second argument in the inner product is an element of ${\cal S}$.
\end{IEEEproof}

\section{Analysis of SLNN}\label{app:slnn}


\begin{IEEEproof}[Solution of the inner subproblem \eqref{eq:costfunction1_2_nn_4_inner}]
For any optimization problem with differentiable objective and constraint functions for which strong duality holds, any set of primal and dual optimal points must satisfy the KKT conditions \cite{Boyd2004}. We define the Lagrangian of \eqref{eq:costfunction1_2_nn_4_inner} with dual variable $\bm{\Lambda}$ as
\begin{equation}\label{eq:costfunction1_2_nn_5}
L(\mathbf{V},\bm{\Lambda}) = \text{tr}(\mathbf{C}^{T}\mathbf{U}\mathbf{V}) + \text{tr}\bigl( \mathbf{\Lambda}^{T}(\mathbf{V}^{T}\mathbf{V}-\mathbf{I}_{n}) \bigr).
\end{equation}
The first-order KKT conditions are given by
\begin{align}
  \nabla_{\mathbf{V}}L(\mathbf{V,\Lambda}) & = \mathbf{U}^{T}\mathbf{C} + \mathbf{V}(\bm{\Lambda}+\bm{\Lambda}^{T}) = \mathbf{0}, \label{eq:costfunction1_2_nn_6}\\
  \mathbf{V}^{T}\mathbf{V} & = \mathbf{I}_{n}, \label{eq:slnn_orthogonality}
\end{align}
where \eqref{eq:costfunction1_2_nn_6} is obtained by setting to zero the gradient\footnote{We use the standard results $\frac{\partial}{\partial\mathbf{X}}\text{tr}(\mathbf{A}^{T}\mathbf{X}) = \mathbf{A}$ and $\frac{\partial}{\partial\mathbf{X}}\text{tr}(\mathbf{X}\mathbf{B}\mathbf{X}^{T}) = \mathbf{X}(\mathbf{B} + \mathbf{B}^{T})$ \cite{Magnus2007,Petersen2008}.} of \eqref{eq:costfunction1_2_nn_5} with respect to $\mathbf{V}$, whereas \eqref{eq:slnn_orthogonality} is the original orthogonality constraint in \eqref{eq:costfunction1_2_nn_4_inner}.

Premultiplying \eqref{eq:costfunction1_2_nn_6} with $\mathbf{V}^{T}$, taking the trace (i.e., taking the inner product with $\mathbf{V}$), and using \eqref{eq:slnn_orthogonality} yields the optimal value for the cost function
\begin{equation}\label{eq:optimal_cost_inner}
  \text{tr}(\mathbf{C}^{T}\mathbf{U}\mathbf{V}) = \text{tr}(\mathbf{V}^{T}\mathbf{U}^{T}\mathbf{C}) = -\text{tr}(\bm{\Lambda}+\bm{\Lambda}^{T}).
\end{equation}
But from $\mathbf{U}^{T}\mathbf{C} = -\mathbf{V}(\bm{\Lambda}+\bm{\Lambda}^{T})$ in \eqref{eq:costfunction1_2_nn_6} we can square both sides to get
\begin{equation}\label{eq:kkt_grad_square}
  \mathbf{C}^{T}\mathbf{U}\mathbf{U}^{T}\mathbf{C} = (\bm{\Lambda}+\bm{\Lambda}^{T})^{2}.
\end{equation}
Hence, among candidate optimal points satisfying the KKT system, the cost function can be made as small as possible by choosing $\bm{\Lambda}+\bm{\Lambda}^{T}$ in \eqref{eq:optimal_cost_inner} as a positive semidefinite matrix square root of the left-hand side of \eqref{eq:kkt_grad_square}. Replacing this in \eqref{eq:optimal_cost_inner} gives the final optimal cost
\begin{equation}
  \text{tr}(\mathbf{C}^{T}\mathbf{U}\mathbf{V}) = -\text{tr}\bigl( (\mathbf{C}^{T}\mathbf{U}\mathbf{U}^{T}\mathbf{C})^{\frac{1}{2}} \bigr) = -\|\mathbf{C}^{T}\mathbf{U}\|_{N}.
\end{equation}
%
\end{IEEEproof}
Interestingly, we point out that the more usual Frobenius norm solves the following relaxed version of the inner subproblem \eqref{eq:costfunction1_2_nn_4_inner}
\begin{equation}\label{eq:costfunction1_2_nn_4_inner_relaxed}
\begin{array}{cl}
\mbox{minimize}& \text{tr}(\mathbf{C}^{T}\mathbf{U}\mathbf{V}) = \langle \mathbf{V}, \mathbf{U}^{T}\mathbf{C} \rangle\\
\mathbf{V}&\\
\mbox{subject to}&\text{tr}(\mathbf{V}^{T}\mathbf{V}) = \|\mathbf{V}\|_{F}^{2} = n,\\
\end{array}
\end{equation}
which is easily verified by writing the KKT system based on the Lagrange function $\text{tr}(\mathbf{C}^{T}\mathbf{U}\mathbf{V}) + \lambda(\text{tr}(\mathbf{V}^{T}\mathbf{V})-n)$,
\begin{align}
  \mathbf{\mathbf{U}^{T}\mathbf{C}} + 2\lambda\mathbf{V} & = 0, & \text{tr}(\mathbf{V}^{T}\mathbf{V}) & = n,
\end{align}
whose solution at the minimum is
\begin{align}
  \mathbf{V} & = -\sqrt{n}\frac{\mathbf{U}^{T}\mathbf{C}}{\|\mathbf{U}^{T}\mathbf{C}\|_{F}}, & \lambda & = \frac{\|\mathbf{U}^{T}\mathbf{C}\|_{F}}{2\sqrt{n}},
\end{align}
with optimal cost $-\sqrt{n}\|\mathbf{U}^{T}\mathbf{C}\|_{F}$. The minimum cost within the expanded domain of this relaxed subproblem will at least be as low as that of \eqref{eq:costfunction1_2_nn_4_inner}, hence $\|\mathbf{U}^{T}\mathbf{C}\|_{N} \leq \sqrt{n}\|\mathbf{U}^{T}\mathbf{C}\|_{F}$. On the other hand,
\begin{equation}
\|\mathbf{U}^{T}\mathbf{C}\|_{N} = \sqrt{\bigl( \sum_{i}\sigma_{i} \bigr)^{2}} \geq \sqrt{\sum_{i}\sigma_{i}^{2}} = \|\mathbf{U}^{T}\mathbf{C}\|_{F},
\end{equation}
where $\sigma_{i}$ denotes the $i$-th singular value of $\mathbf{U}^{T}\mathbf{C}$. Combining the two inequalities we have the bounds
\begin{align}
\|\mathbf{U}^{T}\mathbf{C}\|_{F} \leq \|\mathbf{U}^{T}\mathbf{C}\|_{N} \leq \sqrt{n}\|\mathbf{U}^{T}\mathbf{C}\|_{F}.
\end{align}
\begin{IEEEproof}[Proof of equivalence between \eqref{eq:costfunction1_2_nn_9} and \eqref{eq:costfunction1_2_nn_11}]
We first rewrite \eqref{eq:costfunction1_2_nn_11} replacing the linear matrix inequality with an equivalent Schur complement
\begin{equation}\label{eq:costfunction1_2_nn_11_schur}
\begin{array}{cl}
\mbox{maximize}& 2\,\text{tr}(\mathbf{Z}) + \frac{1}{m}\text{tr}(\mathbf{rr}^T\mathbf{W})\\
\mathbf{W,Z}&\\
\mbox{subject to}&\mathbf{W}\succeq 0, \quad W_{ii}=1\\
& \mathbf{Z}^{2} \preceq \mathbf{C}^{T}\mathbf{W}\mathbf{C}, \quad \mathbf{Z} \succeq 0.
\end{array}
\end{equation}
Let $p_{1}^{*}$ and $p_{2}^{*}$ be the optimal values of problems \eqref{eq:costfunction1_2_nn_9} and \eqref{eq:costfunction1_2_nn_11_schur}, respectively.

Choose a feasible point $(\mathbf{Z},\mathbf{W})$ for the second problem, such that $0 \preceq \mathbf{Z}^2 \preceq \mathbf{C}^{T}\mathbf{W}\mathbf{C}$. This implies\footnote{$\mathbf{A} \succeq \mathbf{B} \succeq 0 \Rightarrow \mathbf{A}^{\frac{1}{2}} \succeq \mathbf{B}^{\frac{1}{2}} \succeq 0$ \cite{Horn1990}.}  $\mathbf{Z} \preceq (\mathbf{C}^{T}\mathbf{W}\mathbf{C})^{\frac{1}{2}}$, hence the values of the two objective functions satisfy
\begin{equation}\label{eq:slnn_sdp_ineq1}
2\,\text{tr}(\mathbf{Z}) + \frac{1}{m}\text{tr}(\mathbf{rr}^{T}\mathbf{W}) \leq 2\,\text{tr}\bigl( (\mathbf{C}^{T}\mathbf{W}\mathbf{C})^{\frac{1}{2}} \bigr) + \frac{1}{m}\text{tr}(\mathbf{rr}^{T}\mathbf{W}).
\end{equation}
In particular, choosing for $(\mathbf{Z},\mathbf{W})$ the unique maximizer of \eqref{eq:costfunction1_2_nn_11_schur}, inequality \eqref{eq:slnn_sdp_ineq1} asserts that $p_{1}^{*} \geq p_{2}^{*}$.

For the converse choose a feasible point $\mathbf{W}$ for the first problem and consider the eigendecomposition $\mathbf{C}^{T}\mathbf{W}\mathbf{C} = \mathbf{Q}\bm{\Lambda}\mathbf{Q}^{T}$. Now set $\mathbf{Z} = \mathbf{Q}\bm{\Lambda}^{\frac{1}{2}}\mathbf{Q}^{T}$, so that $\mathbf{Z}^{2} = \mathbf{Q}\bm{\Lambda}\mathbf{Q}^{T} = \mathbf{C}^{T}\mathbf{W}\mathbf{C}$, and $(\mathbf{W},\mathbf{Z})$ is therefore feasible for \eqref{eq:costfunction1_2_nn_11_schur}. For both problems the value of the cost function is
\begin{equation}
  2\,\text{tr}(\bm{\Lambda}^{\frac{1}{2}}) +  \frac{1}{m}\text{tr}(\mathbf{rr}^{T}\mathbf{W}).
\end{equation}
In particular, choosing for $\mathbf{W}$ the maximizer of \eqref{eq:costfunction1_2_nn_9} the construction for $\mathbf{Z}$ yields a feasible point $(\mathbf{W},\mathbf{Z})$ for \eqref{eq:costfunction1_2_nn_11_schur} where the objective function equals $p_{1}^{*}$. Therefore $p_{1}^{*} \leq p_{2}^{*}$, and coupling this with the converse inequality above we conclude that $p_{1}^{*} = p_{2}^{*}$ and the two problems are equivalent.
\end{IEEEproof}

\bibliographystyle{IEEEtran}
\bibliography{strings1}


\end{document}